\documentclass[12pt]{amsart}
\usepackage{amssymb}
\usepackage[active]{srcltx}
\usepackage{epsfig}
\usepackage{amscd}
\usepackage{a4wide}
\usepackage{amsthm}
\usepackage{tikz}
\usepackage{tikz-cd} 
\usetikzlibrary{matrix}
\newtheorem{theorem}{Theorem}[section]
\newtheorem{lemma}[theorem]{Lemma}
\newtheorem{proposition}[theorem]{Proposition}
\newtheorem{corollary}[theorem]{Corollary}

\theoremstyle{definition}
\newtheorem{definition}[theorem]{Definition}
\newtheorem{example}[theorem]{Example}
\newtheorem{remark}[theorem]{Remark}

\numberwithin{equation}{section}

\DeclareMathOperator{\Sphere}{S}
\newcommand{\bc}{{\mathbb C}}
\newcommand{\bz}{{\mathbb Z}}

\begin{document}
\title [Continuant, Chebyshev polynomials, and Riley polynomials] {Continuant, Chebyshev polynomials, and Riley polynomials}
\author{Kyeonghee Jo and Hyuk Kim}
\address{Division of Liberal Arts and sciences,  Mokpo National Maritime University, Mokpo, Chonnam, 530-729, Korea  }
\email{khjo@mmu.ac.kr}
\address{Department of Mathematical Sciences, Seoul National University, Seoul, 08826, Korea}
\email{hyukkim@snu.ac.kr}

\begin{abstract}
In the previous paper, we showed that the Riley polynomial $\mathcal{R}_K(\lambda)$ of each 2-bridge knot $K$ is split into $\mathcal{R}_K(-u^2)=\pm g(u)g(-u)$, for some integral coefficient polynomial $g(u)\in \mathbb Z[u]$.
In this paper, we study this splitting property of the Riley polynomial.  
 We  show that the Riley polynomial  can be expressed by `$\epsilon$-Chebyshev polynomials', which is a generalization of  Chebyshev polynomials containing the information of $\epsilon_i$-sequence $(\epsilon_i=(-1)^{[i\frac{\beta}{\alpha}]})$ of the 2-bridge knot $K=S(\alpha,\beta)$, and then we give an explicit formula for the splitting polynomial $g(u)$  also as $\epsilon$-Chebyshev polynomials. As applications, we find a sufficient condition for the irreducibility of the Riley polynomials and show the unimodal property of  the symmetrized Riley polynomial.
\end{abstract}

\maketitle

\section{introduction}
Riley investigated non-abelian parabolic representations of 2-bridge knot
groups in \cite{Riley1} and he showed that for each 2-bridge knot $K$, there is an integer coefficient polynomial $\mathcal{R}_K(\lambda) \in \mathbb Z[\lambda]$, the Riley polynomial, and each zero of this polynomial  
corresponds to a parabolic representation
(see Theorem 2 of 
\cite{Riley1}).

In \cite{Jo-Kim},  we studied the  parabolic representations of 2-bridge knots and links by finding arc coloring vectors on the Conway diagram. The method we used is to convert the system of conjugation quandle equations to that of symplectic quandle equations. 
In this approach,  we naturally get to see the Riley polynomial  as $\mathcal{R}_K(-u^2)$ instead of $\mathcal{R}_K(\lambda)$. In the process, we
discovered some interesting arithmetic properties of the Riley polynomials. Namely, splitting of the Riley polynomial, $\mathcal{R}_K(-u^2)=\pm g(u)g(-u)$, for some integral coefficient polynomial $g(u)\in \mathbb Z[u]$, 
 and the fact that the trace field is generated by $u$ instead of $\lambda=-u^2$. (See Theorem 5.1 and Proposition 5.5 of \cite{Jo-Kim}.)

Our main concern in this paper is the splitting property of the Riley polynomial. 
 We will  show in \S\ref{h-formula} that the Riley polynomial  $\mathcal{R}_K(-u^2)$ 
can be expressed by `$\epsilon$-Chebyshev polynomials', which is a generalization of  Chebyshev polynomials containing the information of $\epsilon_i$-sequence $(\epsilon_i=(-1)^{[i\frac{\beta}{\alpha}]})$ of the 2-bridge knot $K=S(\alpha,\beta)$, and then using this together with some properties related to Continuant we'll find an explicit formula for the splitting polynomial $g(u)$ of the Riley polynomial, which is also expressed by $\epsilon$-Chebyshev polynomials. See \S\ref{Cheby} for the definition of $\epsilon$-Chebyshev polynomials.

As applications, we will find a sufficient condition for the irreducibility of the Riley polynomials in \S\ref{irreducibility}, and show  under the assumption of Artin's primitive root conjecture that the set of prime numbers $\alpha$ such that the Riley polynomial of a 2-bridge knot $S(\alpha, \beta)$ is irreducible for each $\beta$ has the asymptotic density $\geq 0.3739558136...$.

In the last section \S\ref{symmetrized Riley}, we will define the symmetrized Riley polynomial for each Riley polynomial and then show that they are all unimodal functions. More precisely, by symmetrizing the variable $u$ to $u=x+\frac{1}{x}$, 
all the Riley polynomials $\mathcal{R}(-u^2)$ become unimodal. This  result is  similar  to the unimodal property of the Alexander polynomial,  well known as the Fox's conjecture.

\section{ Parabolic representations of 2-bridge kmot groups} 
The knot group $G(K)$ of  a 2-bridge knot $K=S(\alpha,\beta)$ always has a presentation of the form 
\begin{equation*}
G(K)=\pi_1(\Sphere^3 \setminus K)=\langle x,y \,|\,wx=yw \rangle,
\end{equation*}
where $w$ is of the form
\begin{equation*}
w=x^{\epsilon_1}y^{\epsilon_2}x^{\epsilon_3}y^{\epsilon_4}\cdots x^{\epsilon_{\alpha-2}}y^{\epsilon_{\alpha-1}}
\end{equation*}
with each $\epsilon_i=(-1)^{[i\frac{\beta}{\alpha}]}$ and  $\epsilon_i=\epsilon_{\alpha-i}$. Here the generators $x$ and $y$ come from the two bridges and represent the meridians.  More generally, Riley defined 2-bridge ``kmot" groups as follows and studied their parabolic representations.
\begin{definition}
A $(2-bridge)$ $kmot$ $group$ is a group $G$ with a presentation  
\begin{equation}\label{kmot-group}
G=\langle x,y \,|\,wx=yw \rangle,
\end{equation}
where $w$ is of the form
\begin{equation}\label{presentation1}
w=x^{\epsilon_1}y^{\epsilon_2}x^{\epsilon_3}y^{\epsilon_4}\cdots x^{\epsilon_{\alpha-2}}y^{\epsilon_{\alpha-1}} \,\,(\alpha \,\text{is an odd integer}\,\geq 3)
\end{equation}
and  $\epsilon_i=\epsilon_{\alpha-i}$ for $i=1, \cdots, \alpha-1$.
\end{definition}
We write ${\bf \epsilon}=(\epsilon_1,\epsilon_2,\cdots,\epsilon_{\alpha-1})$ and $G=G( \epsilon)=G(\epsilon_1,\epsilon_2,\cdots,\epsilon_{\alpha-1})$. 
Suppose that $\rho : G(\epsilon)  \rightarrow SL(2,\bc)$ is a non-abelian $parabolic$ representation, i.e., the trace of $\rho$-image of any meridian is 2. Then after conjugating if necessary, we may assume 
\begin{equation}\label{parabolic-rep}
\rho(x)=\begin{pmatrix}
1& 1\\
0 & 1 
 \end{pmatrix} \quad \text{and}\quad
\rho(y)=\begin{pmatrix}
1& 0\\
\lambda & 1 
 \end{pmatrix}.
\end{equation}
 Riley had shown in \cite{Riley1} that $\lambda$ determines a non-abelian parabolic representation if and only if  $W_{11}=0$. Here $W_{ij}$ is the ($i,j$)-element of 
\begin{equation}\label{parabolic-rep2}
\begin{split}
W(\epsilon)&=\rho(w)=\rho(x)^{\epsilon_1}\rho(y)^{\epsilon_2}\rho(x)^{\epsilon_3}\rho(y)^{\epsilon_4}\cdots \rho(x)^{\epsilon_{\alpha-2}}\rho(y)^{\epsilon_{\alpha-1}}\\
&=\begin{pmatrix}
1& 1\\
0 & 1 
 \end{pmatrix}^{\epsilon_1}\begin{pmatrix}
1& 0\\
\lambda & 1 
 \end{pmatrix}^{\epsilon_2}\begin{pmatrix}
1& 1\\
0 & 1 
 \end{pmatrix}^{\epsilon_3}\cdots \begin{pmatrix}
1& 1\\
0 & 1 
 \end{pmatrix}^{\epsilon_{\alpha-2}}\begin{pmatrix}
1& 0\\
\lambda & 1 
 \end{pmatrix}^{\epsilon_{\alpha-1}},
\end{split}
\end{equation}
for $i,j=1,2$.
Furthermore $W_{11}$ has no repeated roots and the non-abelian parabolic representations bijectively correspond to the roots of the polynomial $W_{11}\in \bz[\lambda]$, the $Riley$ $polynomial$ $\mathcal{R}(\lambda)$. 

Note that  every knot group $G(K)$  of  a 2-bridge knot $K$ is a kmot group, but there is a kmot group which is not a 2-bridge knot group. For example, 
a word $w=xy^{-1}x^{-1}y^{-1}x^{-1}y$ does not define any 2-bridge knot group: For $\alpha=7$, only the following set of $\epsilon_i$-sequences are from 2-bridge knot groups.
$$\{\pm(1,1,1,1,1,1), \pm(1,1,-1,-1,1,1), \pm(1,-1,1,1,-1,1)\}$$

The following formulas are well-known. (See for instance \cite{{M-R}}.)  But we present a proof for the convenience of readers and of later use.
\begin{proposition}\label{W-formula}
For a $\epsilon_i$-sequence ${\bf \epsilon}=(\epsilon_1,\epsilon_2,\cdots,\epsilon_{2m})$, let
\begin{equation*}
\begin{split}
W(\epsilon)&=\begin{pmatrix}
1& 1\\
0 & 1 
 \end{pmatrix}^{\epsilon_1}\begin{pmatrix}
1& 0\\
\lambda & 1 
 \end{pmatrix}^{\epsilon_2}\begin{pmatrix}
1& 1\\
0 & 1 
 \end{pmatrix}^{\epsilon_3}\cdots \begin{pmatrix}
1& 1\\
0 & 1 
 \end{pmatrix}^{\epsilon_{2m-1}}\begin{pmatrix}
1& 0\\
\lambda & 1 
 \end{pmatrix}^{\epsilon_{2m}}=\begin{pmatrix}
{W}_{11} & {W}_{12} \\
{W}_{21}& {W}_{22}
 \end{pmatrix}.
\end{split}
\end{equation*}
Then we have
\begin{equation*}
\begin{split}
W_{11}&= 1+\lambda\sum_{i_1:\text {odd}}^{\wedge} \epsilon_{i_1} \epsilon_{i_2}+ \cdots+\lambda^k\sum_{i_1:\text {odd}}^{\wedge} \epsilon_{i_1} \cdots\epsilon_{i_{2k}}+ \cdots +\lambda^{m} \epsilon_{1} \epsilon_{2}\epsilon_{3}\cdots\epsilon_{2m} \\
W_{12}&= \sum_{i_1:\text {odd}}^{\wedge} \epsilon_{i_1} + \lambda\sum_{i_1:\text {odd}}^{\wedge} \epsilon_{i_1} \epsilon_{i_2}\epsilon_{i_3}+ \cdots+\lambda^{m-1}\epsilon_{1} \epsilon_{2}\epsilon_{3}\cdots\epsilon_{2m-1}   \\
W_{21}&=\lambda\sum_{i_1:\text {even}}^{\wedge} \epsilon_{i_1} + \lambda^2\sum_{i_1:\text {even}}^{\wedge} \epsilon_{i_1} \epsilon_{i_2}\epsilon_{i_3}+ \cdots+\lambda^{m}\epsilon_{2} \epsilon_{3}\epsilon_{4}\cdots\epsilon_{2m}   \\
W_{22}&=1 +\lambda\sum_{i_1:\text {even}}^{\wedge} \epsilon_{i_1} \epsilon_{i_2}+ \lambda^2\sum_{i_1:\text {even}}^{\wedge} \epsilon_{i_1} \epsilon_{i_2}\epsilon_{i_3}\epsilon_{i_4}+ \cdots +\lambda^k\sum_{i_1:\text {even}}^{\wedge} \epsilon_{i_1} \cdots\epsilon_{i_{2k}}+\cdots +\lambda^{m-1} \epsilon_{2} \epsilon_{3} \cdots\epsilon_{2m-1}.
\end{split}
\end{equation*}
Here the notations $\displaystyle\sum_{i_1:\text {odd}}^{\wedge} $ and $\displaystyle\sum_{i_1:\text {even}}^{\wedge} $ are the summations defined as follows.
\begin{equation*}
\begin{split}
\displaystyle\sum_{i_1:\text {odd}}^{\wedge}&  : i_1< i_2 < i_3 <\cdots<i_k, \,\, i_{2i+1} \text{ is odd},\,\,  i_{2i} \text{ is even}\\
\displaystyle\sum_{i_1:\text {even}}^{\wedge}&  : i_1< i_2 < i_3 <\cdots<i_k, \,\, i_{2i+1} \text{ is even},\,\,  i_{2i} \text{ is odd}\\
\end{split}
\end{equation*}
\end{proposition}
\begin{proof}
The formulas hold for $m=1$ since
\begin{equation*}
\begin{split}
W(\epsilon_1, \epsilon_2)&=\begin{pmatrix}
1 & 1\\
  0 & 1
 \end{pmatrix}^{\epsilon_1}\begin{pmatrix}
1 & 0\\
  \lambda& 1 
 \end{pmatrix}^{\epsilon_2}=\begin{pmatrix}
1 & \epsilon_1\\
  0 & 1
 \end{pmatrix}\begin{pmatrix}
1 & 0\\
  \epsilon_2\lambda& 1
 \end{pmatrix}\\
&
=\begin{pmatrix}
1+\lambda \epsilon_1\epsilon_2 & \epsilon_1\\
  \epsilon_2\lambda & 1 
 \end{pmatrix}.
\end{split}
\end{equation*}
Suppose that the formulas  hold for all $m<n$. Then
\begin{equation*}
\begin{split}
W(\epsilon_1, \cdots,\epsilon_{2n})&= W(\epsilon_1, \cdots,\epsilon_{2n-2}) \begin{pmatrix}
1+\lambda \epsilon_{2n-1}\epsilon_{2n} & \epsilon_{2n-1}\\
  \epsilon_{2n}\lambda & 1
 \end{pmatrix},
\end{split}
\end{equation*}
and thus $W(\epsilon_1, \cdots,\epsilon_{2n})_{11}$ is equal to
\begin{equation*}
\begin{split}
&W(\epsilon_1, \cdots,\epsilon_{2n-2})_{11}(1+\lambda \epsilon_{2n-1}\epsilon_{2n} )+W(\epsilon_1, \cdots,\epsilon_{2n-2})_{12} \epsilon_{2n}\lambda \\
&=(1+\lambda \epsilon_{2n-1}\epsilon_{2n} )(1+\lambda\sum_{i_1:\text {odd},\,\, i_2<2n}^{\wedge} \epsilon_{i_1} \epsilon_{i_2}+\cdots+\lambda^{n-1} \epsilon_{1} \epsilon_{2}\cdots\epsilon_{2n-2})\\
&\quad+ \epsilon_{2n}\lambda (\sum_{i_1:\text {odd},\,\, i_1<2n-1}^{\wedge} \epsilon_{i_1} + \lambda\sum_{i_1:\text {odd},\,\, i_3<2n-1}^{\wedge} \epsilon_{i_1} \epsilon_{i_2}\epsilon_{i_3}+ \cdots+\lambda^{n-2}\epsilon_{1} \epsilon_{2}\cdots\epsilon_{2n-3} )\\
&=1 +\lambda (\epsilon_{2n-1}\epsilon_{2n}+\sum_{i_1:\text {odd},\,\, i_2<2n}^{\wedge} \epsilon_{i_1} \epsilon_{i_2}+\sum_{i_1:\text {odd}}^{\wedge} \epsilon_{i_1}\epsilon_{2n} )+ \cdots\\
&\quad+\lambda^{n-1}[(\epsilon_{1} \epsilon_{2}\cdots\epsilon_{2n-2}
 +\epsilon_{2n-1}\epsilon_{2n}(\sum_{i_1:\text {odd},\,\, i_{2n-4}<2n}^{\wedge} \epsilon_{i_1} \cdots \epsilon_{i_{2n-4}})+
\epsilon_{1} \epsilon_{2}\cdots\epsilon_{2n-3}\epsilon_{2n} ]
+\lambda^{n} \epsilon_{1} \epsilon_{2}\cdots\epsilon_{2n}\\
&=1 +\lambda\sum_{i_1:\text {odd}}^{\wedge} \epsilon_{i_1} \epsilon_{i_2}+ \cdots+\lambda^k\sum_{i_1:\text {odd}}^{\wedge} \epsilon_{i_1} \cdots\epsilon_{i_{2k}}+ \cdots +\lambda^{n} \epsilon_{1} \epsilon_{2}\epsilon_{3}\cdots\epsilon_{2n} \\
\end{split}
\end{equation*}
which proves that the formulas hold for $m=n$ also. Similarly we can check that $W(\epsilon_1, \cdots,\epsilon_{2n})_{12}$ and $W(\epsilon_1, \cdots,\epsilon_{2n})_{22}$  also satisfy the given formulas.
\end{proof}
Let $W^*$ be the matrix  obtained from $W$ by exchanging $\begin{pmatrix}
1& 1\\
0 & 1 
 \end{pmatrix}$ and $\begin{pmatrix}
1& 0\\
\lambda & 1 
 \end{pmatrix}$, that is, 
\begin{equation*}
\begin{split}
W^*(\epsilon)&=\begin{pmatrix}
1& 0\\
\lambda & 1 
 \end{pmatrix}^{\epsilon_1}\begin{pmatrix}
1& 1\\
0 & 1 
 \end{pmatrix}^{\epsilon_2}\cdots \begin{pmatrix}
1& 0\\
\lambda & 1 
 \end{pmatrix}^{\epsilon_{2m-1}}\begin{pmatrix}
1& 1\\
0 & 1 
 \end{pmatrix}^{\epsilon_{2m}}.
\end{split}
\end{equation*}
Then we have
\begin{lemma}\label{w-star}
 $ W^*(\epsilon)=\begin{pmatrix}
{W}_{22} & \lambda^{-1}{W}_{21} \\
\lambda {W}_{12}& {W}_{11}
 \end{pmatrix}$.
\end{lemma}
\begin{proof}
The identity holds for $m=1$ since
\begin{equation*}
\begin{split}
W(\epsilon_1, \epsilon_2)&=\begin{pmatrix}
1 & \epsilon_1\\
  0 & 1
 \end{pmatrix}\begin{pmatrix}
1 & 0\\
  \epsilon_2\lambda& 1
 \end{pmatrix}=\begin{pmatrix}
1+\lambda \epsilon_1\epsilon_2 & \epsilon_1\\
  \epsilon_2\lambda & 1
 \end{pmatrix}.
\end{split}
\end{equation*}
and
\begin{equation*}
\begin{split}
W^*(\epsilon_1, \epsilon_2)&=\begin{pmatrix}
1 & 0\\
  \epsilon_1\lambda& 1
 \end{pmatrix}\begin{pmatrix}
1 & \epsilon_2\\
  0 & 1
 \end{pmatrix}=\begin{pmatrix}
1& \epsilon_2\\
  \epsilon_1\lambda & 1 +\lambda \epsilon_1\epsilon_2
 \end{pmatrix}.
\end{split}
\end{equation*}
Assume that the identity  holds for all $m<n$. Then
\begin{equation*}
\begin{split}
W^*(\epsilon_1, \cdots,\epsilon_{2n})&= \begin{pmatrix}
1 & 0\\
  \epsilon_1\lambda & 1
 \end{pmatrix}W(\epsilon_2, \cdots,\epsilon_{2n-1}) \begin{pmatrix}
1 & \epsilon_{2n}\\
  0& 1
 \end{pmatrix}\\
W(\epsilon_1, \cdots,\epsilon_{2n})&= \begin{pmatrix}
1 & \epsilon_1\\
  0 & 1
 \end{pmatrix}W^*(\epsilon_2, \cdots,\epsilon_{2n-1}) \begin{pmatrix}
1 & 0\\
 \epsilon_{2n}\lambda& 1
 \end{pmatrix}.
\end{split}
\end{equation*}
If we let $W'=W(\epsilon_2, \cdots,\epsilon_{2n-1})$, then 
\begin{equation*}
\begin{split}
W^*(\epsilon_1, \cdots,\epsilon_{2n})&= \begin{pmatrix}
1 & 0\\
   \epsilon_1\lambda & 1
 \end{pmatrix}
\begin{pmatrix}
W'_{11} & W'_{12}\\
W'_{21}&  W'_{22}
 \end{pmatrix}
\begin{pmatrix}
1 & \epsilon_{2n}\\
  0& 1
 \end{pmatrix}\\
&=\begin{pmatrix}
W'_{11}& W'_{12} \\
\epsilon_1\lambda W'_{11}+W'_{21}  & \epsilon_1 \lambda W'_{12}+ W'_{22}
 \end{pmatrix}\begin{pmatrix}
1 & \epsilon_{2n}\\
  0& 1
 \end{pmatrix}\\
&=\begin{pmatrix}
W'_{11}& W'_{12}+\epsilon_{2n} W'_{11} \\
\epsilon_1\lambda W'_{11}+W'_{21} & *
 \end{pmatrix}
\end{split}
\end{equation*}
and by the induction hypothesis 
\begin{equation*}
\begin{split}
W(\epsilon_1, \cdots,\epsilon_{2n})&= \begin{pmatrix}
1 & \epsilon_1\\
0   & 1
 \end{pmatrix}\begin{pmatrix}
W'_{22} & \lambda^{-1}W'_{21}\\
\lambda W'_{12}&  W'_{11}
 \end{pmatrix}\begin{pmatrix}
1 & 0\\
  \epsilon_{2n}\lambda& 1
 \end{pmatrix}\\
&=\begin{pmatrix}
W'_{22}+\epsilon_1\lambda  W'_{12} & \lambda^{-1}W'_{21}+\epsilon_1  W'_{11} \\
\lambda W'_{12}& W'_{11}
 \end{pmatrix}\begin{pmatrix}
1 & 0\\
  \epsilon_{2n}\lambda& 1
 \end{pmatrix}\\
&=\begin{pmatrix}
* & \lambda^{-1}W'_{21}+\epsilon_1  W'_{11} \\
\lambda (W'_{12}+\epsilon_{2n} W'_{11})& W'_{11}
 \end{pmatrix}
\end{split}
\end{equation*}
This implies that  the identity  holds for $m=n$ also.
\end{proof}
\begin{corollary}
For a symmetric $\epsilon_i$-sequence ${\bf \epsilon}=(\epsilon_1,\epsilon_2,\cdots,\epsilon_{2m})$, 
that is,  $\epsilon_i=\epsilon_{2m-i+1}$ for all $i=1,2,\cdots,2m$,  we have
\begin{enumerate}
\item [\rm (i)]  $ W(\epsilon)=\begin{pmatrix}
{W}_{11} & {W}_{12} \\
\lambda {W}_{12}& {W}_{22}
 \end{pmatrix}, \,\, W^*(\epsilon)=\begin{pmatrix}
{W}_{22} & {W}_{12} \\
\lambda {W}_{12}& {W}_{11}
 \end{pmatrix}$
\item [\rm (ii)]$\lambda$ defines a non-abelian parabolic representation of the kmot group $G(\epsilon)$ by (\ref{parabolic-rep}) if and only if  $W_{11}=0$.
\end{enumerate}
\end{corollary}
\begin{proof}
(i) follows immediately from Proposition \ref{W-formula}.

Since 
 $$\begin{pmatrix}
{W}_{11} & {W}_{12} \\
\lambda {W}_{12}& {W}_{22}
 \end{pmatrix}\begin{pmatrix}
1& 1\\
0 & 1 
 \end{pmatrix}=
\begin{pmatrix}
{W}_{11} & {W}_{11}+{W}_{12} \\
\lambda {W}_{12}& \lambda {W}_{12}+{W}_{22}
 \end{pmatrix}
$$
and
 $$\begin{pmatrix}
1& 0\\
\lambda & 1 
 \end{pmatrix}\begin{pmatrix}
{W}_{11} & {W}_{12} \\
\lambda {W}_{12}& {W}_{22}
 \end{pmatrix}=\begin{pmatrix}
{W}_{11} & {W}_{12} \\
\lambda ({W}_{11}+{W}_{12})& \lambda {W}_{12}+{W}_{22}
 \end{pmatrix},$$
$$\begin{pmatrix}
{W}_{11} & {W}_{12} \\
\lambda {W}_{12}& {W}_{22}
 \end{pmatrix}\begin{pmatrix}
1& 1\\
0 & 1 
 \end{pmatrix}=\begin{pmatrix}
1& 0\\
\lambda & 1 
 \end{pmatrix}\begin{pmatrix}
{W}_{11} & {W}_{12} \\
\lambda {W}_{12}& {W}_{22}
 \end{pmatrix}$$
is equivalent to $W_{11}=0$, which implies (ii).
\end{proof}

\section{Recursive formula for Riley polynomials}\label{recursive formula}\label{Recursive formula}
For each symmetric $\epsilon_i$-sequence ${\bf \epsilon}=(\epsilon_1,\epsilon_2,\cdots,\epsilon_{2m})$,  we have seen that there exist two monic polynomials $f_{\epsilon}(\lambda)=W_{11}(\lambda),  g_{\epsilon}(\lambda)=W_{12}(\lambda)\in  \mathbb Z[\lambda]$ and $W(\epsilon)\in SL(2,\mathbb Z[\lambda])$ such that 
 $$ W(\epsilon)=\begin{pmatrix}
{W}_{11}(\lambda) & {W}_{12}(\lambda) \\
\lambda {W}_{12}(\lambda)& \displaystyle \frac{1+\lambda{W}_{12}(\lambda)^2}{{W}_{11}(\lambda)}
 \end{pmatrix}=\begin{pmatrix}
f_{\epsilon}(\lambda)& \displaystyle g_{\epsilon}(\lambda) \\
 \lambda g_{\epsilon}(\lambda)& \displaystyle  \frac{1+\lambda g_{\epsilon}(\lambda) ^2}{f_{\epsilon}(\lambda) }
 \end{pmatrix}.$$
Note that 
\begin{equation}\label{f-g}
f_{-\epsilon}(\lambda)=f_{\epsilon}(\lambda), \quad g_{-\epsilon}(\lambda)= -g_{\epsilon}(\lambda). 
\end{equation}

\begin{proposition}\label{recurrences}
Let ${\bf \epsilon}, {\bf \epsilon'}, {\bf \epsilon''}$ be symmetric $\epsilon_i$-sequences such that
$${\bf \epsilon}=( \epsilon_{n}, \cdots, \epsilon_1, \epsilon_1,\cdots, \epsilon_{n}),   {\bf \epsilon'}=( \epsilon_{n-1}, \cdots, \epsilon_1, \epsilon_1, \cdots, \epsilon_{n-1}),  {\bf \epsilon''}=( \epsilon_{n-2}, \cdots,  \epsilon_1, \epsilon_1, \cdots, \epsilon_{n-2}).$$ Then  we have the followings.
\begin{enumerate}
\item [\rm (i)] $g_{\epsilon}= g_{\epsilon'}+\epsilon_n f_{\epsilon'}$
\item [\rm (ii)]
$f_{\epsilon}f_{\epsilon'}=1+\lambda g_{\epsilon}^2$
\item [\rm (iii)] $W=\begin{pmatrix}
f_{\epsilon} & g_{\epsilon}\\
  \lambda g_{\epsilon}& f_{\epsilon'}
 \end{pmatrix}$
\item [\rm (iv)]$f_{\epsilon}=f_{\epsilon''}+\lambda f_{\epsilon'}+ 2\epsilon_n\lambda g_{\epsilon'}$  for $n>1$
\item [\rm (v)]
$\epsilon_n(f_{\epsilon}-f_{\epsilon''})= \lambda (g_{\epsilon}+g_{\epsilon'})$  for $n>1$
\end{enumerate}
Here we let $f_{\epsilon'}(\lambda)=1$ and 
$g_{\epsilon'}(\lambda)=0$
 when ${\bf \epsilon}=(\epsilon_1,\epsilon_1)$.
\end{proposition}
\begin{proof}
(i)-(iv) follow immediately from
\begin{equation*}
\begin{split}
W(\epsilon)&= \begin{pmatrix}
1 & \epsilon_n\\
  0 & 1
 \end{pmatrix}W^*(\epsilon') 
\begin{pmatrix}
1 & 0\\
  \epsilon_n\lambda& 1
 \end{pmatrix}\\
&= \begin{pmatrix}
1 & \epsilon_n\\
  0 & 1
 \end{pmatrix}\begin{pmatrix}
\displaystyle  \frac{1+\lambda g_{\epsilon'}^2}{f_{\epsilon'}} & g_{\epsilon'}\\
  \lambda g_{\epsilon'}& f_{\epsilon'}
 \end{pmatrix}
\begin{pmatrix}
1 & 0\\
 \epsilon_n\lambda& 1
 \end{pmatrix}\\
&= \begin{pmatrix}
\displaystyle  \frac{1+\lambda g_{\epsilon'}^2}{f_{\epsilon'}}+\epsilon_n\lambda g_{\epsilon'}& g_{\epsilon'}+\epsilon_n f_{\epsilon'}\\
 \lambda g_{\epsilon'}& f_{\epsilon'}
 \end{pmatrix}
\begin{pmatrix}
1 & 0\\
  \epsilon_n\lambda& 1
 \end{pmatrix}\\
&= \begin{pmatrix}
\displaystyle  \frac{1+\lambda g_{\epsilon'}^2}{f_{\epsilon'}}+ 2\epsilon_n\lambda g_{\epsilon'}+\lambda f_{\epsilon'}& g_{\epsilon'}+\epsilon_n f_{\epsilon'}\\
 \lambda (g_{\epsilon'}+\epsilon_n f_{\epsilon'})& f_{\epsilon'}
 \end{pmatrix},
\end{split}
\end{equation*}
and  then (v) is proved as follows.
\begin{equation*}
\begin{split}
\epsilon_n f_{\epsilon}&=\epsilon_n(f_{\epsilon''}+\lambda f_{\epsilon'}+ 2\epsilon_n\lambda g_{\epsilon'})\\
&=\epsilon_nf_{\epsilon''}+ \lambda (2g_{\epsilon'}+\epsilon_n f_{\epsilon'})=\epsilon_nf_{\epsilon''}+ \lambda (g_{\epsilon'}+(g_{\epsilon'}+\epsilon_n f_{\epsilon'}))\\
&=\epsilon_nf_{\epsilon''}+ \lambda (g_{\epsilon'}+g_{\epsilon})
\end{split}
\end{equation*}
\end{proof}

Now we can deduce recursive formulas for $f_{ \epsilon}$ and $g_{ \epsilon}$ using Proposition \ref{recurrences} as follows. 
For any sequence $ e_1,  e_2, e_3,\cdots,$ with $ e_i \in\{-1, 1\} $, consider the set of symmetric $ \epsilon_i$-sequences 
$$  {\bf e}(n)=( e_{n}, \cdots, e_2, e_1, e_1, e_2,\cdots, e_{n}),$$
and denote $f_{{\bf e}(n)}$ and $ g_{ {\bf e}(n)}$ by $f_n$ and $g_n$. 
Let $f_0=1$ and $g_0=0$. Then by Proposition \ref{recurrences} we have
\begin{equation}
W( {\bf e}(n))=\begin{pmatrix}
f_{n} & g_{n}\\
  \lambda g_{n} & f_{n-1}
 \end{pmatrix}
\end{equation}
and 
\begin{equation}\label{fNg}
\begin{split}
g_{n}&= g_{n-1}+ e_n f_{n-1}\\
f_{n}f_{n-1}&=1+\lambda g_{n}^2\\
f_{n}&=f_{n-2}+\lambda f_{n-1}+ 2 e_n\lambda g_{n-1}
\end{split}
\end{equation}
Therefore  we have 
\begin{proposition}\label{recursive formulas}
$f_n$ and $g_n$ satisfy the following recursive formulas. 
\begin{enumerate}
\item [\rm (i)] $f_n=f_{n-2}+\lambda f_{n-1}+ 2 e_n\lambda ( e_{n-1} f_{n-2}+ e_{n-2} f_{n-3}+\cdots+ e_2 f_{1}+ e_{1} f_{0})$
\item [\rm (ii)] $g_n= e_n f_{n-1}+ e_{n-1} f_{n-2}+\cdots+ e_2 f_{1}+ e_{1} f_{0}$
\end{enumerate}
\end{proposition}
\begin{proof}
By applying the first identity of (\ref{fNg}) repeatedly, we get (ii): 
\begin{equation*}
\begin{split}
 g_{n}&=g_{n-1}+ e_n f_{n-1}\\
&= e_n f_{n-1}+ e_{n-1} f_{n-2}+g_{n-2}\\
&= e_n f_{n-1}+ e_{n-1} f_{n-2}+\cdots+ e_2 f_{1}+ e_{1} f_{0}
\end{split}
\end{equation*}
and using this and the third identity of (\ref{fNg}) we have (i):
\begin{equation*}
\begin{split}
f_{n}&=f_{n-2}+\lambda f_{n-1}+ 2 e_n\lambda g_{n-1}\\
&=f_{n-2}+\lambda f_{n-1}+ 2 e_n\lambda ( e_{n-1} f_{n-2}+ e_{n-2} f_{n-3}+\cdots+ e_2 f_{1}+ e_{1} f_{0})\\
\end{split}
\end{equation*}
\end{proof}
\section{$\epsilon$-Chebyshev polynomials and the continuant}\label{Cheby}
$Chebyshev$ $polynomials$ are defined by a three-term recursion
$$g_{n+1}(t)=tg_n(t)-g_{n-1}(t).$$ 
We denote the Chebyshev polynomials $g_n(t)$ with the intial condition 
$g_0(t)=a, g_1(t)=b$ by $Ch_n^t(a,b)$, which clearly depends on the initial condition linearly.  And the following properties are also obvious:
\begin{equation}\label{linearity of Cheby}
\begin{split}
Ch_n^t(a,b)&=Ch_{n-k}^t(Ch_k^t(a,b),Ch_{k+1}^t(a,b))\\
Ch_n^t(a,b)&=aCh_n^t(1,0)+bCh_n^t(0,1)
\end{split}
\end{equation}
If we denote a particular Chebyshev polynomials $Ch_n^t(0,1)$ by $s_n(t)$, then 
$$Ch_n^t(1,0)=Ch_{n-1}^t(0,-1)=-Ch_{n-1}^t(0,1)=-s_{n-1}(t),$$
 and arbitrary Chebyshev polynomials are expressed as  linear combinations of $s_{n-1}(t)$ and   $s_n(t)$ as follows :
 $$Ch_n^t(a,b)=-as_{n-1}(t)+bs_n(t)$$
The following identities for Chebyshev polynomials are well-known or easily proved using induction. (See \cite{Hsiao, Rivlin, Yamagishi} for references.)
\begin{lemma}\label{cheby}
Let
\begin{equation*}
\begin{split}
v_n(x)&:=s_{n+1}(x)-s_n(x)=Ch_n^x(1,x-1)=Ch_{n+1}^x(1,1).\\
\end{split}
\end{equation*}
Then
\begin{enumerate}
\item[\rm (i)]   $s_{n}(x)^2-s_{n-1}(x)s_{n+1}(x)=1$
\item [\rm (ii)] $v_n(-x)=(-1)^n (s_{n+1}(x)+s_n(x))$
\item[\rm (iii)] $v_n(x)v_{n-1}(x)=1+(x-2)s_n(x)^2$
\item[\rm (iv)]  $v_{n}(2-x^2)=v_{n}(x)v_{n}(-x)$
\item [\rm (v)] $s_{2n+1}(x)=(-1)^nv_{n}(-x)v_{n}(x)$ 
\item [\rm (vi)]  $\frac{1}{x}s_{2n}(x)=s_n(x^2-2)=(-1)^{n+1}s_n(2-x^2)$
\end{enumerate}
\end{lemma}
Note that $U_n(x)=s_{n+1}(2x)$ and $V_n(x)=v_n(2x)$ are exactly the Chebyshev polynomials of the second and third  kinds, respectively. 
\begin{definition}
For a given $ \epsilon_i$-sequence $\epsilon=(\epsilon_1,\epsilon_2,\epsilon_{3},\cdots), \epsilon_i=\pm1$, 
we will call the  family of polynomials $\{g_n(t)\}$ defined by a three-term recursion
$$g_{n+1}(t)=\epsilon_{n} tg_n(t)-g_{n-1}(t)$$ 
$\epsilon$-$Chebyshev$ $polynomials$. 
\end{definition}
We will denote the $\epsilon$-Chebyshev polynomials $g_n(t)$ with the initial condition 
$g_0(t)=a, g_1(t)=b$ by $\mathcal ECh_n^t(a,b)$, and $\sigma^k(\epsilon)$-Chebyshev polynomials  with the initial condition 
$a, b$ by $\mathcal E^{+k}Ch_n^t(a,b)$, where $\sigma(\epsilon)$ is an $\epsilon_i$-sequence  $1$-shifted from $\epsilon$ to the right, that is,  $$\sigma^k(\epsilon)_i=\epsilon_{i+k}.$$ 
If we denote a particular $\epsilon$-Chebyshev polynomials $\mathcal ECh_n^t(0,1)$ by $s^{\epsilon}_n(t)$, then we have
$$\mathcal ECh_n^t(a,b)=\mathcal E^{+k}Ch_{n-k}^t(\mathcal ECh_k^t(a,b),\mathcal ECh_{k+1}^t(a,b))$$
and
$$\mathcal ECh_n^t(a,b)=-as^{\sigma(\epsilon)}_{n-1}(t)+bs^{\epsilon}_n(t).$$
Note that for a finite  $ \epsilon_i$-sequence $\epsilon=(\epsilon_1,\epsilon_2,\cdots,\epsilon_{m})$, $\mathcal ECh_n^t(a,b)$ is defined only for $n=0,1,\cdots,m+1$ and $\mathcal E^{+k}Ch_n^t(a,b)$ is defined only for $n=0,1,\cdots,m-k+1$. 

The following two lemmas can be proved easily from the definitions of  $\epsilon$-Chebyshev polynomials $s^{\epsilon}_{n}(t)$ and $v^{\epsilon}_{n}(t)$.
\begin{lemma}\label{e-cheby}
Let $v^{\epsilon}_n(t)$ be the $\epsilon$-Chebyshev polynomials defined by
$$
v^{\epsilon}_n(t):=\mathcal ECh_{n+1}^t(1,1).\\
$$
Then
\begin{enumerate}
\item[\rm (i)]  $v^{\epsilon}_{n+2}(t)=\epsilon_{n+2}t v^{\epsilon}_{n+1}(t)-v^{\epsilon}_{n}(t)$, that is, $v^{\epsilon}_n(t)=\mathcal E^{+1}Ch_{n}^t(1,\epsilon_1t-1)$
\item [\rm (ii)]
$v^{\epsilon}_n(t)=s^{\epsilon}_{n+1}(t)-s^{\sigma(\epsilon)}_{n}(t)$
\item [\rm (iii)]  $v^{\epsilon}_n(-t)=(-1)^n(s^{\epsilon}_{n+1}(t)+s^{\sigma(\epsilon)}_{n}(t))=(-1)^n\mathcal ECh_{n+1}^t(-1,1)$
\end{enumerate}
\end{lemma}

\begin{lemma}\label{e-cheby-1}
The following identities hold for $n=0,1,\cdots$.
\begin{enumerate}
\item[\rm (i)]  $s^{-\epsilon}_{n}(t)=(-1)^{n+1}s^{\epsilon}_{n}(t)$
\item [\rm (ii)] $v^{-\epsilon}_{n}(t)=v^{\epsilon}_{n}(-t)$
\end{enumerate}
\end{lemma}

The $continuant$ is a multivariate polynomial representing the determinant of the following tridiagonal matrix. 
$$
K_n(c_1,\cdots,c_n):=\det 
\begin{pmatrix}
c_1 & 1 & & &\\
1&c_2 &1 & &\\
& \ddots & \ddots & \ddots &\\
&& 1&c_{n-1}&1\\
&&&1&c_n
 \end{pmatrix}
$$
We set $K_0:=1, K_{-1}:=0$ for convenience. The continuant has a long history and is related to various fields. For example, continuants occur as both the numerator and the denominator of continued fractions:
$$   
c_1-\frac{1}{c_2-\frac{1}{\ddots\frac{1}{c_{n-1}-\frac{1}{c_n}}}}=\frac{K_n(c_1,\cdots,c_n)}{K_{n-1}(c_2,\cdots,c_n)}
$$
In fact, the term ``continuant"is an abbreviation for "continued-fraction determinant". (See \cite{CO} for a reference.) 

Note that we get the following identities immediately from the definition of the continuant.
\begin{itemize}
\item $K_n(c_1,\cdots,c_n)=c_nK_{n-1}(c_1,\cdots,c_{n-1})-K_{n-2}(c_1,\cdots,c_{n-2})$
\item
$K_n(c_1,\cdots,c_n)=K_n(c_n,\cdots,c_{1})$
\item
$K_n(-c_1,\cdots,-c_n)=(-1)^nK_n(c_1,\cdots,c_n)$
\end{itemize}
The first identity shows the Chebyshev property and we will relate the continuant to $\epsilon$-Chebyshev polynomials by letting $c_i=e_it$.
The following is known as Euler's Identity and play the crucial role for the combinatorial properties of continuant. (See \cite{GKP}, \cite{MO}, or \cite{Ustinov} for references.)
\begin{proposition}[Euler's Identity]\label{Euler Identity}
For $0<i<j<k<l$, the following identity holds.
\begin{equation*}
\begin{split}
&K_{k-i-1}(c_{i+1},\cdots,c_{k-1})K_{l-j-1}(c_{j+1},\cdots,c_{l-1}) \\
=&K_{j-i-1}(c_{i+1},\cdots,c_{j-1})K_{l-k-1}(c_{k+1},\cdots,c_{l-1})+K_{k-j-1}(c_{j+1},\cdots,c_{k-1})K_{l-i-1}(c_{i+1},\cdots,c_{l-1})
\end{split}
\end{equation*}
\end{proposition} 
Applying  Proposition \ref{Euler Identity} to the case $k=j+1$, we get
\begin{corollary}\label{EI-cor}
\begin{equation*}
\begin{split}
&K_{j-i}(c_{i+1},\cdots,c_{j})K_{l-j-1}(c_{j+1},\cdots,c_{l-1})\\
=&K_{j-i-1}(c_{i+1},\cdots,c_{j-1})K_{l-j-2}(c_{j+2},\cdots,c_{l-1})+K_{l-i-1}(c_{i+1},\cdots,c_{l-1})
\end{split}
\end{equation*}
\end{corollary}
Using Corollary \ref{EI-cor} we get the following lemma.
\begin{lemma}\label{EI-1}
The following identities hold.
\begin{enumerate}
\item[\rm (i)] $K_{2n}(c_n, c_{n-1},\cdots,c_1,c_1, c_2, \cdots, c_{n})=K_{n}( c_1,\cdots,c_n)^2-K_{n-1}(c_2,\cdots,c_n)^2$
\item[\rm (ii)] $\sum_{k=1}^n K_{2k}(c_1, c_2, \cdots,c_{k},c_{k}, \cdots, c_2,c_{1})=K_n(c_1, c_2, \cdots, c_{n})^2-1$
\end{enumerate}
\end{lemma}
\begin{proof}
\begin{enumerate}
\item[\rm (i)] Applying  Corollary \ref{EI-cor} to the case $i=0, j=n+1, l=2n+1$, we get
\begin{equation*}
\begin{split}
K_{2n}&(c_n, c_{n-1},\cdots,c_1,c_1, c_2, \cdots, c_{n})\\
&=K_{n}(c_n, \cdots,c_1)K_{n}(c_1, c_{2},\cdots,c_n)-K_{n-1}(c_{n}, \cdots,c_2)K_{n-1}( c_{2},\cdots,c_n)\\
&=K_{n}(c_1, c_{2},\cdots,c_n)^2-K_{n-1}( c_{2},\cdots,c_n)^2
\end{split}
\end{equation*}
\item[\rm (ii)] By applying (i) repeatedly we have
\begin{equation*}
\begin{split}
 K_{n}&(c_1, c_2, \cdots,c_{n})^2
=K_{2n}(c_1, \cdots,c_{n},c_n,\cdots,c_1)+K_{n-1}(c_1, \cdots,c_{n-1})^2\\
&=K_{2n}(c_1, \cdots,c_{n},c_n,\cdots,c_1)+K_{2n-1}(c_1, \cdots,c_{n-1},c_{n-1},\cdots,c_1)+K_{n-2}(c_1, \cdots,c_{n-2})^2\\
&=\sum_{k=n-2}^n K_{2k}(c_1, c_2, \cdots,c_{k},c_{k}, \cdots, c_2,c_{1})+K_{n-3}(c_1, \cdots,c_{n-3})^2\\
&=\cdots\\
&=\sum_{k=2}^n K_{2k}(c_1, c_2, \cdots,c_{k},c_{k}, \cdots, c_2,c_{1})+K_1(c_1)^2\\
&=\sum_{k=1}^n K_{2k}(c_1, c_2, \cdots,c_{k},c_{k}, \cdots, c_2,c_{1})+1\\
\end{split}
\end{equation*}
\end{enumerate}
\end{proof}
\begin{lemma}\label{divisibility1}
The following identities hold.
\begin{enumerate}
\item[\rm (i)]
$K_{2n+1}(c_1,\cdots,c_n, x,-c_n,\cdots,-c_1)=(-1)^nxK_n(c_1,\cdots,c_n)^2$
\item[\rm (ii)]
$\frac{K_{2n+1}(c_1,\cdots,c_n, x,c_1,\cdots,c_n)}{K_n(c_1,\cdots,c_n)}=xK_{n}(c_1,\cdots,c_n)-K_{n-1}(c_2,\cdots,c_{n})-K_{n-1}(c_1,\cdots,c_{n-1})$
\item[\rm (iii)]
$\frac{K_{2n+1}(c_1,\cdots,c_n, x,-c_1,\cdots,-c_n)}{K_n(c_1,\cdots,c_n)}=(-1)^{n}[xK_{n}(c_1,\cdots,c_n)+K_{n-1}(c_2,\cdots,c_{n})-K_{n-1}(c_1,\cdots,c_{n-1})]$
\item[\rm (iv)]
$K_{2n+1}(c_1,\cdots,c_n, x,c_n,\cdots,c_1)=K_n(c_1,\cdots,c_n)(xK_n(c_1,\cdots,c_n)-2K_n(c_1,\cdots,c_{n-1}))$
\end{enumerate}
\end{lemma}
\begin{proof}
Applying  Corollary \ref{EI-cor} to the case $i=0, j=n+1, l=2n+2$, we get the identities (i) and (ii) as follows.
\begin{equation*}
\begin{split}
\rm (i)\,\,&K_{2n+1}(c_1,\cdots,c_n, x,-c_n,\cdots,-c_1)\\
&=K_n(c_1,\cdots,c_n)K_{n+1}(x, -c_n,\cdots,-c_1)-K_{n-1}(c_1,\cdots,c_{n-1})K_{n}( -c_n,\cdots,-c_1)\\
&=(-1)^{n+1}K_n(c_1,\cdots,c_n)K_{n+1}(-x, c_n,\cdots,c_1)-(-1)^nK_{n-1}(c_1,\cdots,c_{n-1})K_{n}(c_n,\cdots,c_1)\\
&=(-1)^{n+1}K_n(c_1,\cdots,c_n)[K_{n+1}(-x, c_n,\cdots,c_1)+K_{n-1}(c_1,\cdots,c_{n-1})]\\
&=(-1)^{n+1}K_n(c_1,\cdots,c_n)[K_{n+1}( c_1,\cdots,c_n,-x)+K_{n-1}(c_1,\cdots,c_{n-1})]\\
&=(-1)^{n+1}K_n(c_1,\cdots,c_n)[-xK_{n}(c_1,\cdots,c_n)-K_{n-1}(c_1,\cdots,c_{n-1})+K_{n-1}(c_1,\cdots,c_{n-1})]\\
&=(-1)^nxK_n(c_1,\cdots,c_n)^2\\
\rm (ii)\,\,&K_{2n+1}(c_1,\cdots,c_n, x,c_1,\cdots,c_n)\\
&=K_n(c_1,\cdots,c_n)K_{n+1}(x, c_1,\cdots,c_n)-K_{n-1}(c_1,\cdots,c_{n-1})K_{n}( c_1,\cdots,c_n)\\
&=K_n(c_1,\cdots,c_n)[K_{n+1}(x, c_1,\cdots,c_n)-K_{n-1}(c_1,\cdots,c_{n-1})]\\
&=K_n(c_1,\cdots,c_n)[xK_{n}(c_n,\cdots,c_1)-K_{n-1}(c_n,\cdots,c_{2})-K_{n-1}(c_1,\cdots,c_{n-1})]\\
&=K_n(c_1,\cdots,c_n)[xK_{n}(c_1,\cdots,c_n)-K_{n-1}(c_2,\cdots,c_{n})-K_{n-1}(c_1,\cdots,c_{n-1})]
\end{split}
\end{equation*}
(iii) and (iv) are also proved similarly.
\end{proof}

Using these identities, we obtain the following divisibility property of continuants.
\begin{corollary}\label{divisibility2}
Let  ${\bf c}=(c_1,\cdots,c_n)$ and $\bar{{\bf c}}=(c_n,\cdots,c_1)$. Then
$$K_n({\bf c})\,|\, K_{nk+n+m}(*,x_1,*,x_2,*,x_3,\cdots,x_m,*)$$
where $*$ is one of $\{{\bf c}, {-\bf c}, \bar{{\bf c}}, -\bar{{\bf c}}\}$.
 \end{corollary}
\begin{proof}
By Lemma \ref{divisibility1} $K_n({\bf c})\,|\,K_{2n+1}(*,x,*)$, 
 and by Corollary \ref{EI-cor} the following identity holds.
\begin{equation*}
\begin{split}
 K&_{nk+n+k}(*,x_1,*,x_2,*,x_3,\cdots,x_k,*)\\
&=K_{n+1}(*,x_1) K_{nk+k-1}(*,x_2,*,x_3,\cdots,x_k,*)-K_n(*)K_{nk+k-2}(*',x_2,*,x_3,\cdots,x_k,*),
\end{split}
\end{equation*}
 where $*'$ is one of 
$\{(c_1,\cdots,c_{n-1}), (c_2,\cdots,c_{n}), (-c_1,\cdots,-c_{n-1}), (-c_2,\cdots,-c_{n})\}$. (For example, if the second $*$ is ${\bf c}=(c_1,\cdots,c_n)$, then $*'=(c_2,\cdots,c_n)$.)  
Hence the statement is easily proved using the induction on $k$.
\end{proof}

The following lemma shows the relationship between the continuant and the  $\epsilon$-Chebyshev polynomials.
\begin{lemma}\label{s-v-continuant}
Let $s^{\epsilon}_n(t)=\mathcal ECh_{n}^t(0,1)$  and $v^{\epsilon}_n(t)=\mathcal ECh_{n+1}^t(1,1)$ be the $\epsilon$-Chebyshev polynomials  for $\epsilon=( e_1, e_2,e_{3},\cdots)$. Then
\begin{enumerate}
\item[\rm (i)]
$s^{\epsilon}_{n+1}(t)=K_{n}( e_1t, e_2t,\cdots,e_nt)$
\item[\rm (ii)]
$v^{\epsilon}_n(t)=K_{n+1}(1, e_1t, e_2t,\cdots,e_nt)=K_{n}( e_1t, e_2t,\cdots,e_nt)-K_{n-1}( e_2t, e_3t,\cdots,e_nt)$
\item[\rm (iii)]
$v^{\epsilon}_n(-t)=(-1)^n[K_{n}( e_1t, e_2t,\cdots,e_nt)+K_{n-1}( e_2t, e_3t,\cdots,e_nt)]$
\end{enumerate}
\end{lemma}
\begin{proof}
(i) follows from the fact that 
$$s^{\epsilon}_1(t)=1=K_0(1), \quad s^{\epsilon}_2(t)=e_1t=K_1(e_1t)$$ 
and 
$$s^{\epsilon}_{k+2}(t)= e_{k+1}ts^{\epsilon}_{k+1}(t)-s^{\epsilon}_k(t),\,\, k\geq 0.$$
Similarly the first equality of (ii) follows from the fact that 
$$v^{\epsilon}_0(t)=1=K_1(1), \quad v^{\epsilon}_1(t)=e_1t-1=K_2(1, e_1t)$$ 
and 
$$v^{\epsilon}_{k+2}(t)= e_{k+2}tv^{\epsilon}_{k+1}(t)-v^{\epsilon}_k(t),\,\, k\geq 0.$$
The second equality of (ii) is obvious.
(iii) is proved as follows.
\begin{equation*}
\begin{split}
v^{\epsilon}_n(-t)&=K_{n}( -e_1t, -e_2t,\cdots,-e_nt)-K_{n-1}( -e_2t, -e_3t,\cdots,-e_nt) \\
&=(-1)^nK_{n}( e_1t, e_2t,\cdots,e_nt)-(-1)^{n-1}K_{n-1}( e_2t, e_3t,\cdots,e_nt)\\
&=(-1)^n[K_{n}( e_1t, e_2t,\cdots,e_nt)+K_{n-1}( e_2t, e_3t,\cdots,e_nt)]
\end{split}    
\end{equation*}
\end{proof}

\begin{corollary}
Let $\epsilon=( e_1, e_2, \cdots,e_{n})$ and $\bar{\epsilon}=( e_n,\cdots,e_2, e_{1})$. Then
\begin{enumerate}
\item[\rm (i)]
$s^{\epsilon}_{n+1}(t)=s_{n+1}^{\bar{\epsilon}}(t)$
\item[\rm (ii)]
$v^{\epsilon}_n(t)=s^{\bar{\epsilon}}_{n+1}(t)-s^{\bar{\epsilon}}_{n}(t)$
\end{enumerate}
\end{corollary}

\begin{corollary}\label{s-v}
Let $\epsilon=(e_1, e_2,\cdots, e_{n})$ and ${\bf e}(n)=( e_{n}, \cdots, e_2, e_1, e_1, e_2,\cdots, e_{n}).$
$$s^{{\bf e}(n)}_{2n+1}(t)=(-1)^nv^{\epsilon}_n(t)v^{\epsilon}_n(-t)$$
\end{corollary}
\begin{proof}
By Lemma \ref{e-cheby}, Lemma \ref{EI-1} and Lemma \ref{s-v-continuant},
\begin{equation*}
\begin{split}
s^{{\bf e}(n)}_{2n+1}(t)&=K_{2n}( e_nt, \cdots, e_1t, e_1t,\cdots,e_{n}t)\\
&=K_{n}( e_1t,\cdots,e_nt)^2-K_{n-1}(e_2t,\cdots,e_nt)^2\\
&=s^{\epsilon}_{n+1}(t)^2-s^{\sigma(\epsilon)}_{n}(t)^2\\
&=(s^{\epsilon}_{n+1}(t)+s^{\sigma(\epsilon)}_{n}(t))(s^{\epsilon}_{n+1}(t)-s^{\sigma(\epsilon)}_{n}(t))\\
&=(-1)^nv^{\epsilon}_n(-t)v^{\epsilon}_n(t)
\end{split}
\end{equation*}
\end{proof}
\section{Riley polynomials as $\epsilon$-Chebyshev polynomials and Divisibility }\label{h-formula}
We have seen in \cite{Jo-Kim} that  the Riley polynomial of any 2-bridge knot $K$ has a special splitting 
 $\mathcal{R}_K(-u^2)=\pm g(u)g(-u)$ for some integral coefficient polynomial $g(u)\in \mathbb Z[u]$, and indeed
 for a Torus knot $K=T(2,2n+1)$,   
\begin{equation}\label{splitting}
{\mathcal R}_K(-u^2)=\pm s_{2n+1}(u)=\pm v_n(u)v_n(-u).
\end{equation}
In this section, we show that (\ref{splitting}) holds true for each 2-bridge knot $S(\alpha,\beta), \alpha=2n+1$ by changing Chebyshev polynomials $v_n(u)$ and $s_{2n+1}(u)$ into $\epsilon$-Chebyshev polynomials as follows.
 $$\mathcal{R}_{S(\alpha,\beta)}(-u^2)=(-1)^{\frac{\alpha-1}{2}}s^{\epsilon}_{\alpha}(u)=(-1)^{n}s^{\epsilon}_{2n+1}(u)=v^{\hat{\epsilon}}_n(u)v^{\hat{\epsilon}}_n(-u)$$
 where
$\epsilon=(\epsilon_1, \epsilon_2, \cdots, \epsilon_{\alpha-1}), \hat{\epsilon}=(\epsilon_{\frac{\alpha-1}{2}}, \epsilon_{\frac{\alpha-1}{2}-1}, \cdots, \epsilon_2, \epsilon_1),\, \epsilon_i=(-1)^{[i\frac{\beta}{\alpha}]}.
 $
 And we also describe a sufficient condition for the divisibility of the Riley polynomials of 2-bridge knots. To do this, we show first that $W_{ij}$ can be expressed as $\epsilon$-Chebyshev polynomials.
\begin{theorem}\label{W-formula2}
Let  ${\bf \epsilon}=(\epsilon_1,\epsilon_2,\cdots,\epsilon_{2n})$, 
 ${\bf \epsilon'}=(\epsilon_2,\cdots,\epsilon_{2n-1})$ and
 \begin{equation*}
\begin{split}
W(\epsilon)&=\begin{pmatrix}
1& 1\\
0 & 1 
 \end{pmatrix}^{\epsilon_1}\begin{pmatrix}
1& 0\\
-u^2 & 1 
 \end{pmatrix}^{\epsilon_2}\begin{pmatrix}
1& 1\\
0 & 1 
 \end{pmatrix}^{\epsilon_3}\cdots \begin{pmatrix}
1& 1\\
0 & 1 
 \end{pmatrix}^{\epsilon_{2n-1}}\begin{pmatrix}
1& 0\\
-u^2 & 1 
 \end{pmatrix}^{\epsilon_{2n}}.
\end{split}
\end{equation*}
Then 
\begin{equation*}
\begin{split}
W(\epsilon)=\begin{pmatrix}
(-1)^ns^{\epsilon}_{2n+1}(u) & \displaystyle (-1)^{n+1}\frac{s^{\epsilon}_{2n}(u)}{u}\\
(-1)^n u  s^{\sigma(\epsilon)}_{2n}(u)& (-1)^{n-1}s^{\epsilon'}_{2n-1}(u)
 \end{pmatrix}.
\end{split}
\end{equation*}
\end{theorem}
\begin{proof}
If we let 
$W(\epsilon)=\begin{pmatrix}
{W}_{11} & {W}_{12} \\
{W}_{21}& {W}_{22}
 \end{pmatrix},$ then by Proposition \ref{W-formula}
\begin{equation*}
\begin{split}
W_{11}&= 1-u^2\sum_{i_1:\text {odd}}^{\wedge} \epsilon_{i_1} \epsilon_{i_2}+ \cdots+(-1)^ku^{2k}\sum_{i_1:\text {odd}}^{\wedge} \epsilon_{i_1} \cdots\epsilon_{i_{2k}}+ \cdots +(-1)^nu^{2n} \epsilon_{1} \epsilon_{2}\epsilon_{3}\cdots\epsilon_{2n} \\
 W_{12}&= \sum_{i_1:\text {odd}}^{\wedge} \epsilon_{i_1} -u^2\sum_{i_1:\text {odd}}^{\wedge} \epsilon_{i_1} \epsilon_{i_2}\epsilon_{i_3}+ \cdots+(-1)^{n-1}u^{2n-2}\epsilon_{1} \epsilon_{2}\epsilon_{3}\cdots\epsilon_{2n-1}   \\
W_{21}&= -u^2\sum_{i_1:\text {even}}^{\wedge} \epsilon_{i_1} +u^4\sum_{i_1:\text {even}}^{\wedge} \epsilon_{i_1} \epsilon_{i_2}\epsilon_{i_3}+ \cdots+(-1)^{n}u^{2n}\epsilon_{2} \epsilon_{3}\epsilon_{4}\cdots\epsilon_{2n}   \\
W_{22}&=1 -u^2\sum_{i_1:\text {even}}^{\wedge} \epsilon_{i_1} \epsilon_{i_2}+ u^4\sum_{i_1:\text {even}}^{\wedge} \epsilon_{i_1} \epsilon_{i_2}\epsilon_{i_3}\epsilon_{i_4}+ \cdots +(-1)^{n-1}u^{2n-2} \epsilon_{2} \epsilon_{3} \cdots\epsilon_{2n-1}.
\end{split}
\end{equation*}
Hence it suffices to show that 
\begin{equation}\label{Eq-11}
W_{11}=(-1)^ns^{\epsilon}_{2n+1}(u)
\end{equation} 
and 
\begin{equation}\label{Eq-12}
W_{12}=(-1)^{n+1}\frac{s^{\epsilon}_{2n}(u)}{u}.
\end{equation}

Firstly, the formulas hold for $n=1$ because
\begin{equation*}
\begin{split}
W(\epsilon_1, \epsilon_2)&=\begin{pmatrix}
1 & 1\\
  0 & 1
 \end{pmatrix}^{\epsilon_1}\begin{pmatrix}
1 & 0\\
 -u^2& 1 
 \end{pmatrix}^{\epsilon_2}=\begin{pmatrix}
1 & \epsilon_1\\
  0 & 1
 \end{pmatrix}\begin{pmatrix}
1 & 0\\
  -\epsilon_2u^2& 1
 \end{pmatrix}\\
&
=\begin{pmatrix}
1-u^2\epsilon_1\epsilon_2 & \epsilon_1\\
  -\epsilon_2u^2& 1 
 \end{pmatrix}
=\begin{pmatrix}
-s^{\epsilon}_3(u) & \frac{s^{\epsilon}_2(u)}{u}\\
  -us^{\sigma(\epsilon)}_2(u)& s^{\epsilon}_1(u) 
 \end{pmatrix}.
\end{split}
\end{equation*}
Suppose that the formulas  hold for all $n<m$ to proceed by induction on $n$. Since
\begin{equation*}
\begin{split}
W(\epsilon_1, \cdots,\epsilon_{2m})&= W(\epsilon_1, \cdots,\epsilon_{2m-2}) 
\begin{pmatrix}
1+\lambda \epsilon_{2m-1}\epsilon_{2m} & \epsilon_{2m-1}\\
  \epsilon_{2m}\lambda & 1
 \end{pmatrix}\\
&=\begin{pmatrix}
(-1)^{m-1}s^{\epsilon}_{2m-1}(u) & \displaystyle (-1)^{m}\frac{s^{\epsilon}_{2m-2}(u)}{u}\\
(-1)^{m-1} u  s^{\sigma(\epsilon)}_{2m-2}(u)& (-1)^{m}s^{\epsilon'}_{2m-3}(u)
 \end{pmatrix}
\begin{pmatrix}
1-u^2\epsilon_{2m-1}\epsilon_{2m} & \epsilon_{2m-1}\\
  -\epsilon_{2m}u^2 & 1
 \end{pmatrix}\\
\end{split}
\end{equation*}
we have the following.
\begin{equation*}
\begin{split}
 W_{11}&=(-1)^m[(u^2\epsilon_{2m-1}\epsilon_{2m}-1)s^{\epsilon}_{2m-1}(u)-\epsilon_{2m}us^{\epsilon}_{2m-2}(u)] \\
 &=(-1)^m[\epsilon_{2m}u(\epsilon_{2m-1}us^{\epsilon}_{2m-1}(u)-s^{\epsilon}_{2m-2}(u))-s^{\epsilon}_{2m-1}(u)]]\\
 &=(-1)^m(\epsilon_{2m}us^{\epsilon}_{2m}(u)-s^{\epsilon}_{2m-1}(u))\\
 &=(-1)^ms^{\epsilon}_{2m+1}(u)\\
 W_{12}&=(-1)^{m-1}\epsilon_{2m-1}s^{\epsilon}_{2m-1}(u)+\displaystyle (-1)^{m}\frac{s^{\epsilon}_{2m-2}(u)}{u}\\
 &=\displaystyle \frac{(-1)^{m-1} }{u}(\epsilon_{2m-1}us^{\epsilon}_{2m-1}(u)-s^{\epsilon}_{2m-2}(u))\\
 &=\displaystyle \frac{(-1)^{m-1} }{u}s^{\epsilon}_{2m}(u)
\end{split}
\end{equation*}
So (\ref{Eq-11}) and  (\ref{Eq-12}) hold  for $n=m$ as well.
\end{proof}

\begin{example}
In the case when all the $\epsilon_i$ are equal to $1$, 
which is the case that the group is the knot group of a Torus knot $K=T(2,2n+1)$,
\begin{equation*}
\begin{split}
W(\epsilon)=\begin{pmatrix}
v_{n}(2+\lambda) & s_{n}(2+\lambda)\\
\lambda s_{n}(2+\lambda)& v_{n-1}(2+\lambda)
 \end{pmatrix}.
\end{split}
\end{equation*}
by the identities
$$(-1)^ns_{2n+1}(u)=v_n(2-u^2)=v_n(2+\lambda)$$ 
and
$$\frac{(-1)^{n+1}}{u}s_{2n}(u)=s_n(2-u^2)=s_n(2+\lambda)$$
 in Lemma \ref{cheby}.
Hence the Riley polynomial of $K$ is equal to 
$$\mathcal R_K(\lambda)=f_n(\lambda)=v_n(2+\lambda).$$
Note that 
\begin{equation}\label{torus-f}
\mathcal R_K(-u^2)=v_n(u)v_n(-u)
\end{equation} 
since $v_n(2+\lambda)=v_n(2-u^2)=v_n(u)v_n(-u)$ 
by Lemma \ref{cheby}, (iv).
\end{example}
Corollary \ref{s-v} and  Theorem \ref{W-formula2} show that (\ref{torus-f}) holds for all the Riley polynomials of kmot groups.
\begin{corollary}\label{cor1}
Let ${\bf e}=(e_1, e_2, \cdots, e_n), \,e_i=\pm 1$ and ${\bf e}(n)=( e_{n}, \cdots, e_2, e_1, e_1, e_2,\cdots, e_{n}).$
Then the Riley polynomial $\mathcal{R}(\lambda)$ of the kmot group $G({\bf e}(n))$ is expressed as 
$$\mathcal{R}(-u^2)=(-1)^ns^{{\bf e}(n)}_{2n+1}(u)=v^{{\bf e}}_n(u)v^{{\bf e}}_n(-u)=f_n(-u^2),$$
where $f_n$ is the polynomial recursively defined in Proposition \ref{recursive formulas}. 
In particular, the  Riley polynomial $\mathcal{R}_K(\lambda)$ of a 2-bridge knot $K=S(\alpha,\beta)$ satisfies  
\begin{equation*}
\begin{split}
\mathcal{R}_K(-u^2)&=(-1)^{\frac{\alpha-1}{2}}s^{\epsilon}_{\alpha}(u),\quad \epsilon=(\epsilon_1, \epsilon_2, \cdots, \epsilon_{\alpha-1}),\,\, \epsilon_i=(-1)^{[i\frac{\beta}{\alpha}]}\\
&=(-1)^{\frac{\alpha-1}{2}}v^{\hat{\epsilon}}_n(u)v^{\hat{\epsilon}}_n(-u),\,\, \hat{\epsilon}=(\epsilon_{\frac{\alpha-1}{2}}, \epsilon_{\frac{\alpha-1}{2}-1}, \cdots, \epsilon_2, \epsilon_1)
\end{split}
\end{equation*}
\end{corollary}
From Corollary \ref{cor1}, we can observe immediately the well-known fact that the Riley polynomial of  any 2-bridge knot is  the same  as its mirror's, since $s^{-{\bf e}(n)}_{2n+1}(u)=s^{{\bf e}(n)}_{2n+1}(u)$ by Lemma \ref{e-cheby-1}. 
Here is another immediate consequence of Corollary \ref{cor1}.  
\begin{corollary}\label{cor2}
The Riley polynomial $\mathcal{R}(\lambda)\in \mathbb Z[\lambda]$ of the kmot group $G({\bf e}(n))$ is irreducible if and only if $v^{{\bf e}}_n(u)\in \mathbb Z[u]$ is irreducible.
\end{corollary}
\begin{proof}
Since $\mathcal{R}(-u^2)=v^{{\bf e}}_n(u)v^{{\bf e}}_n(-u)$ by Corollary \ref{cor1}, 
it is obvious that the reducibility of $v^{{\bf e}}_n(u)$ implies the reducibility of $\mathcal{R}(\lambda)$. 

Now suppose that  $\mathcal{R}(\lambda)=A(\lambda)B(\lambda)$ for some $A(\lambda), B(\lambda)\in \mathbb Z[\lambda]$. Then we have
$$v^{{\bf e}}_n(u)v^{{\bf e}}_n(-u)=\mathcal{R}(-u^2)=A(-u^2)B(-u^2).$$
If $v^{{\bf e}}_n(x)$ is irreducible, then either 
$$
v^{{\bf e}}_n(u)=A(-u^2),\quad v^{{\bf e}}_n(-u) =B(-u^2)$$ 
or
$$
v^{{\bf e}}_n(-u) =A(-u^2),\quad v^{{\bf e}}_n(u) =B(-u^2)$$ 
holds and thus $A(-u^2)=B(-u^2)$. But this implies that $\mathcal{R}(\lambda)=A(\lambda)^2$, which contradicts the fact that the Riley polynomial of any 2-bridge knot has no repeated roots \cite{Riley1}. Therefore $v^{{\bf e}}_n(u)$ must be also reducible and this completes the proof.
\end{proof}

\begin{proposition}\label{div-Riley}
Let $G=G({\bf e})$ and $\tilde{G}=G(\tilde{{\bf e}})$ be two kmot groups. Then the Riley polynomial  of $\tilde{G}$ has that of $G$ as a factor if $\tilde{{\bf e}}$ is expressed as
\begin{equation}\label{5-4}
\tilde{{\bf e}}=({\bf e}, \delta_1, \pm {\bf e}, \delta_2, \pm {\bf e},\delta_3,\cdots,\pm {\bf e}, \delta_{2k},{\bf e})    
\end{equation}
where $\delta_i\in\{1,-1\}$ for $i=1,2,\cdots,2k$.
\end{proposition}
\begin{proof}
Let ${\bf e}=(e_1,\cdots,e_{2n})$ and ${\bf c}=(e_1u,\cdots,e_{2n}u)$. Then the Riley polynomials $\mathcal{R}_G$ and $\mathcal{R}_{\tilde{G}}$ of $G$ and $\tilde{G}$ satisfy the following identities by Lemma \ref{s-v-continuant} and Corollary \ref{cor1}. (Note that $\bar{{\bf e}}={\bf e}$ and  $\bar{\tilde{{\bf e}}}=\tilde{{\bf e}}$ since both $G$ and $\tilde{G}$ are kmot groups.)
\begin{equation*}
\begin{split}
\mathcal{R}_G(-u^2)&=(-1)^ns^{{\bf e}}_{2n+1}(u)=(-1)^nK_{2n}(c_1,\cdots,c_{2n})=(-1)^nK_{2n}({\bf c})\\
\mathcal{R}_{\tilde{G}}(-u^2)&=(-1)^{(2k+1)n+k}s^{\tilde{{\bf e}}}_{(2k+1)(2n+1)}(u)\\
&=(-1)^{(2k+1)n+k}K_{2k(2n+1)+2n}({\bf c},\delta_1u,\pm {\bf c},\delta_2u, \pm {\bf c},\cdots,\pm {\bf c},\delta_{2k}u,{\bf c})\\
\end{split}
\end{equation*}
Hence $\mathcal{R}_{\tilde{G}}(\lambda)$ has $\mathcal{R}_G(\lambda)$ as a factor by Corollary \ref{divisibility2}.
\end{proof}
\begin{example}
For the $\epsilon_i$-sequence ${\bf e}=(1,1)$ of the trefoil $3_1=S(3,1)$, there are 4 $\epsilon_i$-sequences with length $8$ which satisfy the identity (\ref{5-4}):
\begin{equation*}
\begin{split}
(1,1,1,1,1,1,1,1),(1,1,-1,1,1,-1,1,1), (1,1,-1,-1,-1,-1,1,1),(1,1,1,-1,-1,1,1,1)
\end{split}
\end{equation*}
The first one $(1,1,1,1,1,1,1,1)$ is the $\epsilon_i$-sequence of 2-bridge knot $9_1=S(9,1)$, and the remaining three in the second line are all $\epsilon_i$-sequences of kmot groups which are not from 2-bridge knots.
\end{example}

\begin{remark}
Note that if $G_1=G$ and $G_2=\tilde G$ in Proposition \ref{div-Riley} are the knot groups of 2-bridge knots $K_1$ and $K_2$ respectively, then the divisibility of Riley polynomials implies that there exists an epimorphism from $G_2$ to $G_1$ by Kitano-Morifuji \cite{KM}. Furthermore,  if  $C[a_1, a_2,\cdots,a_m]$ is the Conway normal form of $K_1$ and ${\bf a}=(a_1, a_2,\cdots,a_m)$, then by Aimi-Lee-Sakai-Sakuma \cite{ALS} $K_2$ has an ORS-expansion  with respect to ${\bf a}$, that is,
$$K_2=C[\delta_1{\bf a}, 2c_1, \delta_2 {\bf a^{-1}}, 2c_2, \delta_3 {\bf a},2c_3, \delta_4 {\bf a^{-1}},2c_4,\cdots, \delta_{2n} {\bf a^{-1}},2c_{2n}, \delta_{2n+1} {\bf a}]$$ 
where
$${\bf a^{-1}}=(a_m, a_{m-1},\cdots,a_1),\ \delta_i=\pm 1\ (\delta_1=1), \ c_i \in \mathbb Z.$$

The converse statement of Proposition \ref{div-Riley} is not true in general. That is, 
the existence of an epimorphism  between two knot groups does not imply that the corresponding $\epsilon_i$-sequences of the two 2-bridge knots have such a relation as in Proposition  \ref{div-Riley}. 
For example,  there exists an epimorphism from the knot group of $9_6$ to that of $3_1=S(3,1)$,
but the $\epsilon_i$-sequences of $S(27,5)$ and $S(27,11)$ are 
$$(1,1,1,1,1,-1,-1,-1,-1,-1,1,1,1,1,1,1,-1,-1,-1,-1,-1,1,1,1,1,1)$$
and
$$(1,1,-1,-1,1,1,1,-1,-1,1,1,1,-1,-1,1,1,1,-1,-1,1,1,1,-1,-1,1,1)$$
respectively. Note that $S(27,5)$ is the Schubert form of $9_6$ and $S(27,11)$ is that of its upside-down in the Conway normal form expression. 
\end{remark}

Using Lemma \ref{EI-1}, Lemma \ref{s-v-continuant}, and Corollary \ref{cor1}, we have another recursive formula for the Riley polynomials of kmot groups as follows.
\begin{corollary}
Let ${\bf e'}(n)$ be a symmetric sequence $  {\bf e'}(n)=( e_1, \cdots, e_{n-1}, e_{n}, e_{n}, e_{n-1},\cdots, e_1).$ Then the sequence of the Riley polynomials $\mathcal{R}_k(\lambda)$ of  the kmot group $G({\bf e'}(k))$ satisfies
\begin{equation*}
\begin{split}
\sum_{k=1}^n (-1)^k\mathcal{R}_k(-u^2)=K_n(e_1u, \cdots, e_{n-1}u, e_{n}u)^2-1=s^{\epsilon}_{n+1}(u)^2-1
\end{split}
\end{equation*} 
where $ \epsilon=(e_1, e_2, \cdots, e_n)$.
\end{corollary}
\begin{proof}
\begin{equation*}
\begin{split}
s^{\epsilon}_{n+1}(u)^2-1&=K_n(e_1u, e_2u, \cdots, e_nu)^2-1\\
&=\sum_{k=1}^n K_{2k}(e_1u, e_2u, \cdots, e_ku, e_ku,  \cdots, e_{2}u,, e_1u)\\
&=\sum_{k=1}^n (-1)^k\mathcal{R}_k(-u^2)
\end{split}
\end{equation*} 
\end{proof}
\section{Irreducibility of Riley polynomials}\label{irreducibility}

Let $\Phi_k(x)$ be the $k$-th cyclotomic polynomial. Then there is a unique polynomial $\Psi_k(x)\in \mathbb Z[x]$ for $k\geq 3$ such that 
$$\Psi_k(x+\frac{1}{x})=x^{-\frac{1}{2}\varphi(k)}\Phi_k(x),$$ 
where $\varphi$ is the Euler function. M. Yamagishi has shown in \cite{Yamagishi} that 
\begin{equation}\label{Yamagishi}
  v_n(x)=\prod_{1<d\mid 2n+1} \Psi_{2d}(x).   
\end{equation}
(See Proposition 2.4, (iii) in \cite{Yamagishi}.) This identity 
implies that the Riley polynomial of the torus knot $T(2,2n+1)=S(2n+1,1)$ is irreducible over $\mathbb Q$ if and only if $2n+1$ is a prime number. 
For arbitrary kmot groups with $\epsilon_i$-sequence of length $2n$, we can prove that their Riley polynomials are all irreducible if $2n+1$ is a prime number and $2$ is a primitive root modulo $2n+1$. 
\begin{theorem}\label{primitive}
Let $\alpha$ be a prime number and $\epsilon=(\epsilon_1,\epsilon_2,\cdots,\epsilon_{\alpha-1}) $ be any symmetric $\epsilon_i$-sequence, that is,   
$\epsilon_i=\epsilon_{\alpha-i}\in\{1,-1\}$ for $i=1, \cdots, \alpha-1$.
Then the Riley polynomial of a kmot group $G(\epsilon)$ is irreducible if $2$ is a primitive root modulo $\alpha$.
\end{theorem}
\begin{proof}
Let  $\alpha=2n+1$  and ${\bf e}=(\epsilon_{n+1},\cdots,\epsilon_{2n})=(\epsilon_{n},\cdots,\epsilon_2,\epsilon_1)$.

By Corollary \ref{cor1} and Corollary \ref{cor2},  the Riley polynomial $\mathcal{R}(\lambda)$ of $G({\bf e}(n))$ satisfies
$$\mathcal{R}(-u^2)=(-1)^ns^{\epsilon}_{2n+1}(u)=v^{{\bf e}}_n(u)v^{{\bf e}}_n(-u),$$
 and
 $\mathcal{R}(\lambda)\in \mathbb Z[\lambda]$  is irreducible if and only if $v^{{\bf e}}_n(u)\in \mathbb Z[u]$ is irreducible.

Hence it suffices to show that  $v_n(u)\in \mathbb Z_2[u]$ is irreducible if $2$ is a primitive root modulo $2n+1$ by the following 2 facts.
\begin{itemize}
\item $v^{{\bf e}}_n(u)$ is equal to $v_n(u)$  in $\mathbb Z_2[u]$. 
\item  $v^{{\bf e}}_n(u)$ is irreducible in $\mathbb Z[u]$ if  $v^{{\bf e}}_n(u)$ is irreducible in $\mathbb Z_2[u]$. 
\end{itemize}

Firstly, we claim that  $v_n(u)$ is equal to $\Psi_{2n+1}(u)$ in $\mathbb Z_2[u]$. 
Since $2n+1$ is a prime, $$ v_n(u)=\prod_{1<d\mid 2n+1} \Psi_{2d}(u)=\Psi_{2(2n+1)}(u)$$ 
by (\ref{Yamagishi}) and thus we have 
\begin{equation}\label{6-2}
v_n(u)=\Psi_{2(2n+1)}(u)=(-1)^{\frac{1}{2}\phi(2n+1)}\Psi_{2n+1}(-u)\equiv \Psi_{2n+1}(u)\in \mathbb Z_2[u],
\end{equation}
which proves our claim. The second equality of (\ref{6-2}) follows from the property
$$
\Phi_{2(2n+1)}(x)=\Phi_{2n+1}(-x).
$$

Since $\Psi_{2n+1}(u)\in \mathbb Z_2[u]$ is irreducible if $\Phi_{2n+1}(x)$ is irreducible in $\mathbb Z_2[x]$,
and
$\Phi_{2n+1}(x)$ is irreducible in $\mathbb Z_2[x]$ if and only if $2$ is a primitive root modulo $2n+1$ as is well known, 
it is immediate that
 $v_n(u)\in \mathbb Z_2[u]$ is irreducible if $2$ is a primitive root modulo $2n+1$, which completes the proof.
\end{proof}
As an immediate consequence of Theorem \ref{primitive} we get the following, which 
Morifuji and Tran have already proved in a different way  in   \cite{M-Tran}.
\begin{corollary}[Morifuji-Tran \cite{M-Tran}]\label{S2}
Let $\alpha$ be a prime number and $\beta$ be any odd integer. Then the Riley polynomial of a 2-bridge knot $S(\alpha, \beta)$ is irreducible if $2$ is a primitive root modulo $\alpha$.
\end{corollary}
Denote by $S(2)$ the set of prime numbers $p$ such that $2$ is a primitive root modulo $p$.  Then for any $\alpha \in S(2)$, the Riley polynomial of a 2-bridge knot $S(\alpha, \beta)$ is irreducible for each $\beta$ by Corollary \ref{S2}.

 Artin's primitive root conjecture claims that the set $S(2)$ is infinite and has the asymptotic density $0.3739558136...$, which is called Artin constant \cite{Moree}. 
 From this, we have the following.
\begin{corollary}
Under the assumption of Artin's conjecture, 
the set of prime numbers $\alpha$ such that the Riley polynomial of a 2-bridge knot $S(\alpha, \beta)$ is irreducible for each $\beta$ has the asymptotic density $\geq 0.3739558136...$.
\end{corollary}

For each  symmetric $\epsilon_i$-sequence $\epsilon=(\epsilon_1,\epsilon_2,\cdots,\epsilon_{\alpha-1}) $, denote the number of irreducible factors of the Riley polynomial of a kmot group $G(\epsilon)=G(\epsilon_1,\epsilon_2,\cdots,\epsilon_{\alpha-1})$ by $n(\epsilon)=n(\epsilon_1,\epsilon_2,\cdots,\epsilon_{\alpha-1})$. Then 
from the above argument  in the proof of Theorem \ref{primitive}, we have the following.
\begin{proposition}
Let  $\alpha$ be an odd number such that $2$ is a primitive root modulo $d\neq 1$ for any  divisor $d$ of $\alpha$.
Then 
$$n(\epsilon)\leq |\{d \in \mathbb N \,|\, d\neq  1, d\mid \alpha\}|$$
for any  symmetric $\epsilon_i$-sequence $\epsilon=(\epsilon_1,\epsilon_2,\cdots,\epsilon_{\alpha-1}) $ of length $\alpha-1$. For the case when $\epsilon_i=1 $ for all $i=1,\cdots,\alpha-1$, that is, for the torus knot $T(2,\alpha)$,  the equality holds.
\end{proposition}
\begin{proof}
Let  $\alpha=2n+1$  and ${\bf e}=(\epsilon_{n+1},\cdots,\epsilon_{2n})=(\epsilon_{n},\cdots,\epsilon_2,\epsilon_1)$. Then $n(\epsilon)$ is equal to the number of irreducible factors of $v^{\bf e}_{n}$ by Corollary \ref{cor1}.
Since $2$ is a primitive root modulo $d\neq 1$ for any  divisor $d$ of $\alpha$, $\Psi_{2d}(x)\equiv \Psi_{d}(x)$ is irreducible in $\mathbb Z_2[x]$ and 
$v_n(x)=\prod_{1<d\mid 2n+1} \Psi_{2d}(x)$, the number of irreducible factors of $v_{n}(x)\equiv v^{\bf e}_{n}(x)$  in $\mathbb Z_2[x]$ is equal to $n(\epsilon)$.
So the number of irreducible factors of $v^{\bf e}_{n}$ in $\mathbb Z[x]$ is less than or equal to $n(\epsilon)$.
\end{proof}
\begin{example}
Let $\alpha=9$. Then $2$ is a primitive root modulo $3$ and $9$, and
$$ |\{d \in \mathbb N \,|\, d\neq  1, d\mid 9\}|=|\{3,9\}|=2.$$
If we denote $v^{\epsilon}_4$  for an $\epsilon_i$-sequence  $\epsilon=(\epsilon_1,\epsilon_2,\epsilon_3,\epsilon_4)$ by $v^{\epsilon_1,\epsilon_2,\epsilon_3,\epsilon_4}_4$, then we have 
$$v^{\epsilon_1,\epsilon_2,\epsilon_3,\epsilon_4}_4(u)=\epsilon_4\epsilon_3\epsilon_2\epsilon_1u^4-\epsilon_4\epsilon_3\epsilon_2u^3-(\epsilon_4\epsilon_3+\epsilon_4\epsilon_1+\epsilon_2\epsilon_1)u^2+(\epsilon_4+\epsilon_2)u+1.$$ 
Hence we can check that $n(\epsilon)\leq 2$ for all $16$ $\epsilon_i$-sequences of length 4 as follows.
\begin{equation*}
\begin{split}
v^{1,1,1,1}_4(u)&=v_4(u)=(u-1)(u^3-3u-1)=\Psi_3(-u)\Psi_9(-u)=\Psi_6(u)\Psi_{18}(u)\\
&=v^{-1,-1,-1,-1}_4(-u)\\
v^{1,1,1,-1}_4(u)&=-(u+1)(u^3-2u^2+u-1)=v^{-1,-1,-1,1}_4(-u)\\
v^{1,1,-1,1}_4(u)&=-u^4+u^3-u^2+2u+1=v^{-1,-1,1,-1}_4(-u)\\
v^{1,-1,1,1}_4(u)&=-(u-1)(u^3+u+1)=v^{-1,1,-1,-1}_4(-u)\\
v^{-1,1,1,1}_4(u)&=-(u+1)(u^3-u-1)=v^{1,-1,-1,-1}_4(-u)\\
v^{-1,-1,1,1}_4(u)&=(u+1)(u^3-u+1)=v^{1,1,-1,-1}_4(-u)\\
v^{1,-1,-1,1}_4(u)&=u^4-u^3+u^2+1=v^{-1,1,1,-1}_4(-u)\\
v^{-1,1,-1,1}_4(u)&=u^4+u^3+3u^2+2u+1=v^{1,-1,1,-1}_4(-u)\\
\end{split}
\end{equation*}
Note that the last two equality imply that the Riley polynomials of 2-bridge knots $S(9,5)$ and  $S(9,7)$ are irreducible.
\end{example}

\section{symmetrized Riley polynomials}\label{symmetrized Riley}
Putting the two generators of a 2-bridege knot group $G(K)=\langle x,y \,|\,wx=yw \rangle$ 
\begin{equation}\label{nonabel-rep}
x=\begin{pmatrix}
M & 1\\
  0 & \frac{1}{M} 
 \end{pmatrix} \quad \text{and}\quad
y=\begin{pmatrix}
M & 0\\
  \lambda& \frac{1}{M} 
 \end{pmatrix},
\end{equation}
all the non-abelian $SL(2,\bc)$ representations are given by the solutions of the 2-variable Riley polynomial  
$$\mathcal{R}(\lambda,s)=0, s=M^2+\frac{1}{M^2}-2$$ 
in $\mathbb Z[\lambda,s]$. (See \cite{Riley2}.) Note that $s=0$ gives the parabolic representations, that is, $\mathcal{R}(\lambda,0)$ is equal to the Riley polynomial $\mathcal{R}(\lambda)$.

In the case of $\lambda=0$, if we symmetrize the variable $s$, more precisely 
$$s+2=t+\frac{1}{t}, t=M^2,$$
then we get the Alexander polynomial $\bigtriangleup_K(t)$ and $\bigtriangleup_K(-t)$ has been proved to be unimodal by  Hartley in \cite{ Hartley}, as is well known for the Fox's conjecture. Here, the unimodality of a polynomial  $\sum_{i=0}^n a_i^n x^i \in \mathbb R[x]$ means that there exists an index $0\leq j\leq n$ for which one of the following two properties holds:
\begin{enumerate}
\item[\rm (i)] $a_0\leq a_1\leq a_2\leq\cdots\leq a_{j-1}\leq a_j\geq a_{j+1}\geq a_{j+2}\geq\cdots\geq a_n $
\item[\rm (ii)] $a_0\geq a_1\geq a_2\geq\cdots\geq a_{j-1}\geq a_j\leq a_{j+1}\leq a_{j+2}\leq\cdots\leq a_n $
\end{enumerate} 
Furthermore if the sequence $\{a_0, a_1, \cdots, a_n\}$ has 
the strict monotonicity, we will say $\sum_{i=0}^n a_i^n x^i$ is $strongly$ $unimodal$. 

Interestingly, we can see that  the same result can be achieved  by reversing the roles of the two variables $s$ and $\lambda=-u^2$. By symmetrizing the variable $u$ to $u=x+\frac{1}{x}$ when $s=0$,
the Riley polynomial $\mathcal{R}(-u^2)$ becomes to be unimodal.

Since  $s^{\epsilon}_{n+1}(u)=K_{n}( \epsilon_1u, \epsilon_2u,\cdots,\epsilon_nu)$ is either an odd polynomial or an even polynomial for any $\epsilon=( \epsilon_1, \epsilon_2,\cdots)$, we have a palindromic polynomial $\tilde{\phi}^{\epsilon}_n(x)$ such that
$$\tilde{\phi}^{\epsilon}_n(x^2)=x^{n}s^{\epsilon}_{n+1}(x+\frac{1}{x}).$$ 
Note that
$$\tilde{\phi}^{\epsilon}_0(x)=1,\,\,\tilde{\phi}^{\epsilon}_1(x)=\epsilon_1(x+1)$$
and the following identity follows from the recursion of $s^{\epsilon}_{n+1}(u)$.
$$\tilde{\phi}^{\epsilon}_{n+1}(x)=\epsilon_{n+1}(x+1)\tilde{\phi}^{\epsilon}_n(x)-x\tilde{\phi}^{\epsilon}_{n-1}(x)$$
If we let 
$$\phi^{\epsilon}_{n}(x):=\epsilon_1\epsilon_2\cdots \epsilon_n\tilde{\phi}^{\epsilon}_{n}(x),$$
then we have 
$$\phi^{\epsilon}_0(x)=1,\,\,\phi^{\epsilon}_1(x)=x+1$$
and
\begin{equation}\label{phi-identity}
\phi^{\epsilon}_{n+1}(x)=(x+1)\phi^{\epsilon}_n(x)-\epsilon_{n}\epsilon_{n+1}x\phi^{\epsilon}_{n-1}(x).
\end{equation}
Now we prove that $\phi^{\epsilon}_{n}$ is  unimodal, which implies that  $\tilde{\phi}^{\epsilon}_n$ is also unimodal. Note that if $n$ is even and $\epsilon$ is symmetric, then $\phi^{\epsilon}_{n}=\tilde{\phi}^{\epsilon}_{n}$. 
\begin{proposition}\label{symRiley-prop}
Let $\phi^{\epsilon}_n(x)=\sum_{j=0}^n a_j^n x^j$. Then 
\begin{enumerate}
\item[\rm (i)] $a_0^n=a_n^n=1$, $a_j^n\geq 0$ for  $j=0,1,\cdots,n$.
\item[\rm (ii)] $a_j^n\geq a_{j-1}^{n}$  for  $j=0,1,\cdots,[\frac{n}{2}]$.
\item[\rm (iii)]  $a_j^n-a_{j-1}^{n}\geq a_j^{n-1}-a_{j-1}^{n-1}$  for $j=0,1,\cdots,[\frac{n}{2}]$.
\item[\rm (iv)] $a_j^n\geq a_j^{n-1}$ for $j=0,1,\cdots,n-1$.
\item[\rm (v)] $a_{1}^n=a_{n-1}^n=1+2s$, where $s$ is the number of sign changes in  $(\epsilon_1, \epsilon_2,\cdots, \epsilon_n)$.
\item[\rm (vi)] $ a_{k}^n=a_{n-k}^n=\sum_{m=0}^k(-1)^m\begin{pmatrix}n-2m\\k-m  \end{pmatrix}\sum_{1\leq j_i<j_{i+1}-1<n-1}\delta_{j_1}\cdots\delta_{j_m}$ for $n>1$ and $k=0,1,\cdots,[\frac{n}{2}]$. Here $\delta_i=\epsilon_i\epsilon_{i+1}$.
\end{enumerate}
\end{proposition}
\begin{proof}
From (\ref{phi-identity}),
$$\sum_{j=0}^{n+1} a_j^{n+1} x^j=(x+1)\sum_{j=0}^n a_j^n x^j-\epsilon_{n}\epsilon_{n+1}x\sum_{j=0}^{n-1} a_j^{n-1} x^j$$
and thus
\begin{equation}\label{coeff-identity}
a_j^{n+1}=a_j^{n}+a_{j-1}^{n}-\epsilon_{n}\epsilon_{n+1}a_{j-1}^{n-1},\,\,j=1,2,\cdots, n.
\end{equation}
Hence  we have
$$a_j^{n+1}-a_{j-1}^{n+1}=(a_j^{n}-a_{j-1}^{n})+(a_{j-1}^{n}-a_{j-2}^{n})-\epsilon_{n}\epsilon_{n+1}(a_{j-1}^{n-1}-a_{j-2}^{n-1})$$ for $j=1,2,\cdots, n$.

Since $\phi^{\epsilon}_2(x)=(x+1)\phi^{\epsilon}_1(x)-\epsilon_2\epsilon_1x=(x+1)^2-\epsilon_2\epsilon_1x=x^2+(2-\epsilon_2\epsilon_1)x+1$,
the statements (i)-(vi) hold for $n=1, 2$. Suppose that they hold for all $n\leq k$ to proceed by induction on $n$. Then  for  $j=0,1,\cdots,[\frac{k+1}{2}]$ 
\begin{equation*}
\begin{split}
a_j^{k+1}-a_{j-1}^{k+1}&= 
(a_j^{k}-a_{j-1}^{k})+(a_{j-1}^{k}-a_{j-2}^{k})-\epsilon_{k}\epsilon_{k+1}(a_{j-1}^{k-1}-a_{j-2}^{k-1})\\
&\geq (a_j^{k}-a_{j-1}^{k})+(a_{j-1}^{k}-a_{j-2}^{k})-(a_{j-1}^{k-1}-a_{j-2}^{k-1})\\
&\geq a_j^{k}-a_{j-1}^{k}\\
&\geq 0
\end{split}
\end{equation*}
and  for $j=0,1,\cdots,k$
\begin{equation*}
\begin{split}
a_j^{k+1}-a_j^{k}&=a_{j-1}^{k}-\epsilon_{k}\epsilon_{k+1}a_{j-1}^{k-1} \quad \\
&\geq a_{j-1}^{k}-a_{j-1}^{k-1} \quad \\
&\geq 0
\end{split}
\end{equation*}
which imply the statements (i)-(iv) hold for $n=k+1$ as well.

To prove (v), it suffices to show that 
\begin{equation}\label{2nd coeff}
a_{n-1}^n=n-\sum_{j=1}^{n-1}\epsilon_i\epsilon_{i+1}\,\,\, \text{for}\,\,n>1
\end{equation}
because
$$n-\sum_{j=1}^{n-1}\epsilon_i\epsilon_{i+1}=1+\sum_{j=1}^{n-1}(1-\epsilon_i\epsilon_{i+1})=1+2s.$$
Since $\phi^{\epsilon}_2(x)=x^2+(2-\epsilon_2\epsilon_1)x+1$, $a_1^2=2- \epsilon_2\epsilon_1$ and thus (\ref{2nd coeff}) holds for $n=2$.
If we assume that (\ref{2nd coeff}) holds for all $n\leq k$, then 
\begin{equation*}
\begin{split}
a_k^{k+1}&=a_k^{k}+a_{k-1}^{k}-\epsilon_{k}\epsilon_{k+1}a_{k-1}^{k-1}\\
&=1+(k-\sum_{j=1}^{k-1}\epsilon_i\epsilon_{i+1})-\epsilon_{k}\epsilon_{k+1}\\
&=(k+1)-\sum_{j=1}^{k}\epsilon_i\epsilon_{i+1},
\end{split}
\end{equation*}
which implies that  (\ref{2nd coeff}) holds for $n=k+1$ as well.

By (\ref{coeff-identity}) and (\ref{2nd coeff}), we have 
$$a_j^{n+1}=a_j^{n}+a_{j-1}^{n}-\delta_{n}a_{j-1}^{n-1},\,\,j=1,2,\cdots, n$$
and
\begin{equation*}
a_{n-1}^n=n-\sum_{j=1}^{n-1}\epsilon_i\epsilon_{i+1}=\begin{pmatrix}n\\1  \end{pmatrix}-\begin{pmatrix}n\\0  \end{pmatrix}\sum_{j=1}^{n-1}\delta_i.
\end{equation*}
Using these two identities  and the identity 
$$\begin{pmatrix}n+1-2m\\k-m  \end{pmatrix}=\begin{pmatrix}n-2m\\k-m  \end{pmatrix}+\begin{pmatrix}n-2m\\k-m-1 \end{pmatrix}
$$
(vi) is proved by the induction argument as follows.
\begin{equation*}
\begin{split}
\sum_{m=0}^k&(-1)^m\begin{pmatrix}{n+1}-2m\\k-m  \end{pmatrix}\sum_{1\leq j_i<j_{i+1}-1<n}\delta_{j_1}\cdots\delta_{j_m}\\
&=\sum_{m=0}^k(-1)^m(\begin{pmatrix}n-2m\\k-m  \end{pmatrix}+\begin{pmatrix}n-2m\\k-m-1 \end{pmatrix})\sum_{1\leq j_i<j_{i+1}-1<n}\delta_{j_1}\cdots\delta_{j_m}\\
&=\sum_{m=0}^k(-1)^m\begin{pmatrix}n-2m\\k-m  \end{pmatrix}\sum_{1\leq j_i<j_{i+1}-1<n-1}\delta_{j_1}\cdots\delta_{j_m}\\
&\quad+\sum_{m=0}^k(-1)^m\begin{pmatrix}n-2m\\(k-1)-m  \end{pmatrix}\sum_{1\leq j_i<j_{i+1}-1<n-1}\delta_{j_1}\cdots\delta_{j_m}\\
&\quad+\delta_n\sum_{m=0}^k(-1)^m\begin{pmatrix}{n+1}-2m\\k-m  \end{pmatrix}\sum_{1\leq j_i<j_{i+1}-1<n-2}\delta_{j_1}\cdots\delta_{j_{m-1}}\\
&= a_{n-k}^{n}+ a_{n-(k-1)}^{n}-\delta_n\sum_{m=0}^k(-1)^{m-1}\begin{pmatrix}{n+1}-2m\\k-m  \end{pmatrix}\sum_{1\leq j_i<j_{i+1}-1<n-2}\delta_{j_1}\cdots\delta_{j_{m-1}}\\
&= a_{n-k}^{n}+ a_{n-(k-1)}^{n}-\delta_n\sum_{m=0}^k(-1)^{m-1}\begin{pmatrix}(n-1)-2(m-1)\\(k-1)-(m-1)  \end{pmatrix}\sum_{1\leq j_i<j_{i+1}-1<n-2}\delta_{j_1}\cdots\delta_{j_{m-1}}\\
&= a_{k-1}^{n}+ a_{(n+1)-k}^{n}-\delta_na_{n-k}^{n-1}(= a_{k}^{n}+ a_{k-1}^{n}-\delta_na_{k-1}^{n-1} \,\text{by the induction hypothesis})\\
&= a_{(n+1)-k}^{n+1} (= a_{k}^{n+1})
\end{split}
\end{equation*}
\end{proof}
\begin{corollary}
Let $\epsilon=( \epsilon_1, \epsilon_2,\cdots, \epsilon_{n}),\,\, \epsilon_i=\pm 1$.  Then
\begin{enumerate}
\item[\rm (i)] 
The coefficients of $\phi_n^{\epsilon}(x)$ are minimum when $\epsilon_1=\epsilon_i$  for all $i=1,2,\cdots,n$. In this case, 
$\phi_n^{\epsilon}(x)=\sum_{j=0}^n  x^j.$
\item[\rm (ii)] 
The coefficients of $\phi_n^{\epsilon}(x)$ are maximum  when $\epsilon_i$ is alternating.
\item[\rm (iii)] $\phi_n^{\epsilon}(x)$ is stronly unimodal if $\epsilon_i$ is alternating.
\end{enumerate}
\end{corollary}
\begin{proof}
If  $\epsilon_i=\epsilon_1$  for all $i=1,2,\cdots,n$, then
$$\phi_0(x)=1,\,\,\phi_1(x)=x+1,\,\,\phi_{n+1}(x)=(x+1)\phi_n(x)-x\phi_{n-1},$$ 
so (i) is proved by the induction argument.
And (ii) follows from the following inequality.
\begin{equation*}
\begin{split}
a_{n-k}^n&=\sum_{m=0}^k(-1)^m\begin{pmatrix}n-2m\\k-m  \end{pmatrix}\sum_{1\leq j_i<j_{i+1}-1<n-1
}\delta_{j_1}\cdots\delta_{j_m}\\
&\leq \sum_{m=0}^k\begin{pmatrix}n-2m\\k-m  \end{pmatrix}\sum_{{\tiny  \begin{array}{cc}
( j_1,\cdots,j_m)\\1\leq j_i<j_{i+1}-1<n-1
\end{array}}}1\\
&=\sum_{m=0}^k(-1)^m\begin{pmatrix}n-2m\\k-m  \end{pmatrix}\sum_{{\tiny  \begin{array}{cc}
( j_1,\cdots,j_m)\\1\leq j_i<j_{i+1}-1<n-1
\end{array}}}(-1)^m \quad (\text{when}\,\,\delta_i=-1\,\, \text{for all }\,\, i=1,\cdots,n-1)\\
\end{split}
\end{equation*}
(iii) is obvious from the identity 
$$a_{k}^n=a_{n-k}^n= \sum_{m=0}^k\begin{pmatrix}n-2m\\k-m  \end{pmatrix}\sum_{{\tiny  \begin{array}{cc}
( j_1,\cdots,j_m)\\1\leq j_i<j_{i+1}-1<n-1
\end{array}}}1,\quad k=0,1,\cdots,[\frac{n}{2}].$$
\end{proof}

\begin{definition}
When $\epsilon=( \epsilon_1, \epsilon_2,\cdots, \epsilon_{2n})$ is the $\epsilon_i$-sequence of a 2-bridge kmot group $G(\epsilon)$, the Riley polynomial of $G(\epsilon)$ is $\mathcal{R}(-u^2)=(-1)^ns_{2n+1}^{\epsilon}(u)=(-1)^nK_{2n}( \epsilon_1u, \epsilon_2u,\cdots,\epsilon_{2n}u)$. And 
we will call the   polynomial $\phi^{\epsilon}_{2n}(x)$ the {\it symmetrized Riley polynomial} of $G(\epsilon)$, where
$\phi^{\epsilon}_{2n}(x^2)=x^{2n}s^{\epsilon}_{2n+1}(x+\frac{1}{x}).$
\end{definition}
\begin{corollary}
The symmetrized Riley polynomial of any 2-bridge kmot group is unimodal.
\end{corollary}

This symmetrized Riley polynomial is actually obtained from $\mathcal{R}(\lambda, s)$ by substituting 
\begin{equation*}\label{sym-1}
s=0,\,\, -(\lambda+2)=X+\frac{1}{X},\,\, X=x^2
\end{equation*} similar to the Alexander polynomial $\bigtriangleup_K(-t)$  obtained by substituting 
\begin{equation*}\label{sym-2}
\lambda=0,\,\, -(s+2)=T+\frac{1}{T}, \,\,T=-t.
\end{equation*}
That is, symmetrizing 2-variable Riley polynomial  $\mathcal{R}(\lambda, s)$ by 
$$-(\lambda+2)=X+\frac{1}{X},\,\, -(s+2)=T+\frac{1}{T}, $$
we get a polynomial $F(T,X)$ such that $$F(T,-1)=\bigtriangleup_K(-t),\,\,F(-1,X)=\phi^{\epsilon}_{\alpha-1}(X).$$
Note that $\lambda+2$ is the trace of $\rho(x)\rho(y)=\begin{pmatrix}
1& 1\\
0 & 1 
 \end{pmatrix}\begin{pmatrix}
1& 0\\
\lambda & 1 
 \end{pmatrix}$ and $s+2$ is the trace of $\rho(x)^2=\begin{pmatrix}
M& 1\\
0 & 1/M 
 \end{pmatrix}^2$ as well as $\rho(y)^2=\begin{pmatrix}
M& 0\\
\lambda & 1/M 
 \end{pmatrix}^2$. 

\section*{Acknowledgements }
Both authors were supported 
 by the National Research Foundation of Korea(NRF) grant funded by the Korea government(MSIT) (NRF-2021R1F1A1049444 and NRF-2018R1A2B6005691, respectively). The authors would like to 
thank the anonymous referee for his helpful comments and Philip Choi for his  help with computer calculations.

\end{document}